\expandafter\chardef\csname pre amssym.def at\endcsname=\the\catcode`\@ 
\catcode`\@=11 
 
\def\undefine#1{\let#1\undefined} 
\def\newsymbol#1#2#3#4#5{\let\next@\relax 
 \ifnum#2=\@ne\let\next@\msafam@\else 
 \ifnum#2=\tw@\let\next@\msbfam@\fi\fi 
 \mathchardef#1="#3\next@#4#5} 
\def\mathhexbox@#1#2#3{\relax 
 \ifmmode\mathpalette{}{\m@th\mathchar"#1#2#3}%
 \else\leavevmode\hbox{$\m@th\mathchar"#1#2#3$}\fi} 
\def\hexnumber@#1{\ifcase#1 0\or 1\or 2\or 3\or 4\or 5\or 6\or 7\or 8\or 
 9\or A\or B\or C\or D\or E\or F\fi} 
 
\font\tenmsa=msam10 
\font\sevenmsa=msam7 
\font\fivemsa=msam5 
\newfam\msafam 
\textfont\msafam=\tenmsa 
\scriptfont\msafam=\sevenmsa 
\scriptscriptfont\msafam=\fivemsa 
\edef\msafam@{\hexnumber@\msafam} 
\mathchardef\dabar@"0\msafam@39 
\def\dashrightarrow{\mathrel{\dabar@\dabar@\mathchar"0\msafam@4B}} 
\def\dashleftarrow{\mathrel{\mathchar"0\msafam@4C\dabar@\dabar@}} 
 
\def\ulcorner{\delimiter"4\msafam@70\msafam@70 } 
\def\urcorner{\delimiter"5\msafam@71\msafam@71 } 
\def\llcorner{\delimiter"4\msafam@78\msafam@78 } 
\def\lrcorner{\delimiter"5\msafam@79\msafam@79 } 
\def\yen{{\mathhexbox@\msafam@55 }} 
\def\checkmark{{\mathhexbox@\msafam@58 }} 
\def\circledR{{\mathhexbox@\msafam@72 }} 
\def\maltese{{\mathhexbox@\msafam@7A }} 
 
\font\tenmsb=msbm10 
\font\sevenmsb=msbm7 
\font\fivemsb=msbm5 
\newfam\msbfam 
\textfont\msbfam=\tenmsb 
\scriptfont\msbfam=\sevenmsb 
\scriptscriptfont\msbfam=\fivemsb 
\edef\msbfam@{\hexnumber@\msbfam}

\catcode`\@=\csname pre amssym.def at\endcsname 
 
\expandafter\ifx\csname pre amssym.tex at\endcsname\relax \else \endinput\fi 
\expandafter\chardef\csname pre amssym.tex at\endcsname=\the\catcode`\@ 
\catcode`\@=11 
\newsymbol\boxdot 1200 
\newsymbol\boxplus 1201 
\newsymbol\boxtimes 1202 
\newsymbol\square 1003 
\newsymbol\blacksquare 1004 
\newsymbol\centerdot 1205 
\newsymbol\lozenge 1006 
\newsymbol\blacklozenge 1007 
\newsymbol\circlearrowright 1308 
\newsymbol\circlearrowleft 1309 
\undefine\rightleftharpoons 
\newsymbol\rightleftharpoons 130A 
\newsymbol\leftrightharpoons 130B 
\newsymbol\boxminus 120C 
\newsymbol\Vdash 130D 
\newsymbol\Vvdash 130E 
\newsymbol\vDash 130F 
\newsymbol\twoheadrightarrow 1310 
\newsymbol\twoheadleftarrow 1311 
\newsymbol\leftleftarrows 1312 
\newsymbol\rightrightarrows 1313 
\newsymbol\upuparrows 1314 
\newsymbol\downdownarrows 1315 
\newsymbol\upharpoonright 1316 
  
\newsymbol\downharpoonright 1317 
\newsymbol\upharpoonleft 1318 
\newsymbol\downharpoonleft 1319 
\newsymbol\rightarrowtail 131A 
\newsymbol\leftarrowtail 131B 
\newsymbol\leftrightarrows 131C 
\newsymbol\rightleftarrows 131D 
\newsymbol\Lsh 131E 
\newsymbol\Rsh 131F 
\newsymbol\rightsquigarrow 1320 
\newsymbol\leftrightsquigarrow 1321 
\newsymbol\looparrowleft 1322 
\newsymbol\looparrowright 1323 
\newsymbol\circeq 1324 
\newsymbol\succsim 1325 
\newsymbol\gtrsim 1326 
\newsymbol\gtrapprox 1327 
\newsymbol\multimap 1328 
\newsymbol\therefore 1329 
\newsymbol\because 132A 
\newsymbol\doteqdot 132B 
  
\newsymbol\triangleq 132C 
\newsymbol\precsim 132D 
\newsymbol\lesssim 132E 
\newsymbol\lessapprox 132F 
\newsymbol\eqslantless 1330 
\newsymbol\eqslantgtr 1331 
\newsymbol\curlyeqprec 1332 
\newsymbol\curlyeqsucc 1333 
\newsymbol\preccurlyeq 1334 
\newsymbol\leqq 1335 
\newsymbol\leqslant 1336 
\newsymbol\lessgtr 1337 
\newsymbol\backprime 1038 
\newsymbol\risingdotseq 133A 
\newsymbol\fallingdotseq 133B 
\newsymbol\succcurlyeq 133C 
\newsymbol\geqq 133D 
\newsymbol\geqslant 133E 
\newsymbol\gtrless 133F 
\newsymbol\sqsubset 1340 
\newsymbol\sqsupset 1341 
\newsymbol\vartriangleright 1342 
\newsymbol\vartriangleleft 1343 
\newsymbol\trianglerighteq 1344 
\newsymbol\trianglelefteq 1345 
\newsymbol\bigstar 1046 
\newsymbol\between 1347 
\newsymbol\blacktriangledown 1048 
\newsymbol\blacktriangleright 1349 
\newsymbol\blacktriangleleft 134A 
\newsymbol\vartriangle 134D 
\newsymbol\blacktriangle 104E 
\newsymbol\triangledown 104F 
\newsymbol\eqcirc 1350 
\newsymbol\lesseqgtr 1351 
\newsymbol\gtreqless 1352 
\newsymbol\lesseqqgtr 1353 
\newsymbol\gtreqqless 1354 
\newsymbol\Rrightarrow 1356 
\newsymbol\Lleftarrow 1357 
\newsymbol\veebar 1259 
\newsymbol\barwedge 125A 
\newsymbol\doublebarwedge 125B 
\undefine\angle 
\newsymbol\angle 105C 
\newsymbol\measuredangle 105D 
\newsymbol\sphericalangle 105E 
\newsymbol\varpropto 135F 
\newsymbol\smallsmile 1360 
\newsymbol\smallfrown 1361 
\newsymbol\Subset 1362 
\newsymbol\Supset 1363 
\newsymbol\Cup 1264 
  
\newsymbol\Cap 1265 
  
\newsymbol\curlywedge 1266 
\newsymbol\curlyvee 1267 
\newsymbol\leftthreetimes 1268 
\newsymbol\rightthreetimes 1269 
\newsymbol\subseteqq 136A 
\newsymbol\supseteqq 136B 
\newsymbol\bumpeq 136C 
\newsymbol\Bumpeq 136D 
\newsymbol\lll 136E 
  
\newsymbol\ggg 136F 
  
\newsymbol\circledS 1073 
\newsymbol\pitchfork 1374 
\newsymbol\dotplus 1275 
\newsymbol\backsim 1376 
\newsymbol\backsimeq 1377 
\newsymbol\complement 107B 
\newsymbol\intercal 127C 
\newsymbol\circledcirc 127D 
\newsymbol\circledast 127E 
\newsymbol\circleddash 127F 
\newsymbol\lvertneqq 2300 
\newsymbol\gvertneqq 2301 
\newsymbol\nleq 2302 
\newsymbol\ngeq 2303 
\newsymbol\nless 2304 
\newsymbol\ngtr 2305 
\newsymbol\nprec 2306 
\newsymbol\nsucc 2307 
\newsymbol\lneqq 2308 
\newsymbol\gneqq 2309 
\newsymbol\nleqslant 230A 
\newsymbol\ngeqslant 230B 
\newsymbol\lneq 230C 
\newsymbol\gneq 230D 
\newsymbol\npreceq 230E 
\newsymbol\nsucceq 230F 
\newsymbol\precnsim 2310 
\newsymbol\succnsim 2311 
\newsymbol\lnsim 2312 
\newsymbol\gnsim 2313 
\newsymbol\nleqq 2314 
\newsymbol\ngeqq 2315 
\newsymbol\precneqq 2316 
\newsymbol\succneqq 2317 
\newsymbol\precnapprox 2318 
\newsymbol\succnapprox 2319 
\newsymbol\lnapprox 231A 
\newsymbol\gnapprox 231B 
\newsymbol\nsim 231C 
\newsymbol\ncong 231D 
\newsymbol\diagup 231E 
\newsymbol\diagdown 231F 
\newsymbol\varsubsetneq 2320 
\newsymbol\varsupsetneq 2321 
\newsymbol\nsubseteqq 2322 
\newsymbol\nsupseteqq 2323 
\newsymbol\subsetneqq 2324 
\newsymbol\supsetneqq 2325 
\newsymbol\varsubsetneqq 2326 
\newsymbol\varsupsetneqq 2327 
\newsymbol\subsetneq 2328 
\newsymbol\supsetneq 2329 
\newsymbol\nsubseteq 232A 
\newsymbol\nsupseteq 232B 
\newsymbol\nparallel 232C 
\newsymbol\nmid 232D 
\newsymbol\nshortmid 232E 
\newsymbol\nshortparallel 232F 
\newsymbol\nvdash 2330 
\newsymbol\nVdash 2331 
\newsymbol\nvDash 2332 
\newsymbol\nVDash 2333 
\newsymbol\ntrianglerighteq 2334 
\newsymbol\ntrianglelefteq 2335 
\newsymbol\ntriangleleft 2336 
\newsymbol\ntriangleright 2337 
\newsymbol\nleftarrow 2338 
\newsymbol\nrightarrow 2339 
\newsymbol\nLeftarrow 233A 
\newsymbol\nRightarrow 233B 
\newsymbol\nLeftrightarrow 233C 
\newsymbol\nleftrightarrow 233D 
\newsymbol\divideontimes 223E 
\newsymbol\varnothing 203F 
\newsymbol\nexists 2040 
\newsymbol\Finv 2060 
\newsymbol\Game 2061 
\newsymbol\mho 2066 
\newsymbol\eth 2067 
\newsymbol\eqsim 2368 
\newsymbol\beth 2069 
\newsymbol\gimel 206A 
\newsymbol\daleth 206B 
\newsymbol\lessdot 236C 
\newsymbol\gtrdot 236D 
\newsymbol\ltimes 226E 
\newsymbol\rtimes 226F 
\newsymbol\shortmid 2370 
\newsymbol\shortparallel 2371 
\newsymbol\smallsetminus 2272 
\newsymbol\thicksim 2373 
\newsymbol\thickapprox 2374 
\newsymbol\approxeq 2375 
\newsymbol\succapprox 2376 
\newsymbol\precapprox 2377 
\newsymbol\curvearrowleft 2378 
\newsymbol\curvearrowright 2379 
\newsymbol\digamma 207A 
\newsymbol\varkappa 207B 
\newsymbol\Bbbk 207C 
\newsymbol\hslash 207D 
\undefine\hbar 
\newsymbol\hbar 207E 
\newsymbol\backepsilon 237F 
\catcode`\@=\csname pre amssym.tex at\endcsname 
 
\magnification=1200 
\hsize=468truept 
\vsize=646truept 
\voffset=-10pt 
\parskip=1pc 
\baselineskip=14truept 
\count0=1 
 
\dimen100=\hsize 
 
\def\leftill#1#2#3#4{ 
\medskip 
\line{$ 
\vcenter{ 
\hsize = #1truept \hrule\hbox{\vrule\hbox to  \hsize{\hss \vbox{\vskip#2truept 
\hbox{{\copy100 \the\count105}: #3}\vskip2truept}\hss } 
\vrule}\hrule} 
\dimen110=\dimen100 
\advance\dimen110 by -36truept 
\advance\dimen110 by -#1truept 
\hss \vcenter{\hsize = \dimen110 
\medskip 
\noindent { #4\par\medskip}}$} 
\advance\count105 by 1 
} 
\def\rightill#1#2#3#4{ 
\medskip 
\line{ 
\dimen110=\dimen100 
\advance\dimen110 by -36truept 
\advance\dimen110 by -#1truept 
$\vcenter{\hsize = \dimen110 
\medskip 
\noindent { #4\par\medskip}} 
\hss \vcenter{ 
\hsize = #1truept \hrule\hbox{\vrule\hbox to  \hsize{\hss \vbox{\vskip#2truept 
\hbox{{\copy100 \the\count105}: #3}\vskip2truept}\hss } 
\vrule}\hrule} 
$} 
\advance\count105 by 1 
} 
\def\midill#1#2#3{\medskip 
\line{$\hss 
\vcenter{ 
\hsize = #1truept \hrule\hbox{\vrule\hbox to  \hsize{\hss \vbox{\vskip#2truept 
\hbox{{\copy100 \the\count105}: #3}\vskip2truept}\hss } 
\vrule}\hrule} 
\dimen110=\dimen100 
\advance\dimen110 by -36truept 
\advance\dimen110 by -#1truept 
\hss $} 
\advance\count105 by 1 
} 
\def\insectnum{\copy110\the\count120 
\advance\count120 by 1 
}

\font\ninerm=cmr9 
\font\eightrm=cmr8

\font\tenrm=cmr10 at 10pt 
 
\font\sc=cmcsc10

 
\def\msb{\fam\msbfam\tenmsb} 
 
\def\bba{{\msb A}} 
 
\def\bbc{{\msb C}}

\def\bbi{{\msb I}}

\def\bbp{{\msb P}} 
\def\bbq{{\msb Q}} 
\def\bbr{{\msb R}}

\def\bbz{{\msb Z}}

\def\grD{\Delta}

\def\grL{\Lambda} 
\def\grO{\Omega}

\def\gra{\alpha}

\def\grd{\delta} 
\def\gre{\epsilon}

\def\gri{\iota} 
 
\def\grl{\lambda}

\def\gro{\omega}

\def\grr{\rho}

\def\grz{\zeta} 
 
\def\la#1{\hbox to #1pc{\leftarrowfill}} 
\def\ra#1{\hbox to #1pc{\rightarrowfill}} 
 
\def\fract#1#2{\raise4pt\hbox{$ #1 \atop #2 $}} 
\def\decdnar#1{\phantom{\hbox{$\scriptstyle{#1}$}} 
\left\downarrow\vbox{\vskip15pt\hbox{$\scriptstyle{#1}$}}\right.} 
 
\def\bowtie{\hbox to 1pt{\hss}\raise.66pt\hbox{$\scriptstyle{>}$} 
\kern-4.9pt\triangleleft} 
\def\hsmash{\triangleright\kern-4.4pt\raise.66pt\hbox{$\scriptstyle{<}$}} 
\def\boxit#1{\vbox{\hrule\hbox{\vrule\kern3pt 
\vbox{\kern3pt#1\kern3pt}\kern3pt\vrule}\hrule}}

\def\za{\vrule height6pt width4pt depth1pt}

\font\aa=eufm10

\def\Got#1{\hbox{\aa#1}}

\def\bfw{{\bf w}}

\def\bfz{{\bf z}}

\def\calc{{\cal C}} 
\def\calo{{\cal O}} 
 
\def\cald{{\cal D}} 
 
\def\calf{{\cal F}}

\def\cali{{\cal I}}

\def\call{{\cal L}} 
\def\calm{{\cal M}} 
 
\def\calo{{\cal O}}

\def\cals{{\cal S}}

\def\calz{{\cal Z}}

\def\gF{{\Got F}} 
\def\gG{{\Got G}}

\def\Got#1{\hbox{\aa#1}}

\def\gsp1{{\Got s}{\Got p}(1)}

\font\svtnrm=cmr17

\font\bsc=cmcsc10 at 10truept

\def\ker{\hbox{ker}}
\def\im{\hbox{im}}
\def\Ric{\hbox{Ric}}
\def\barj{\bar{j}}
\def\Se{Sasakian-Einstein }

\def\intro{0}
\def\trans{1}
\def\alg{4}
\def\link{3}
\def\top{6}
\def\est{2}
\def\pr{5}
\def\mod{7}
\phantom{ooo}
\bigskip\bigskip
\centerline{\svtnrm On the Geometry of \Se}
\bigskip
\centerline{\svtnrm 5-Manifolds}

\bigskip\bigskip
\centerline{\sc Charles P. Boyer~~ Krzysztof Galicki~~ Michael Nakamaye}
\footnote{}{\ninerm During the preparation of this work the first two authors 
were partially supported by NSF grant DMS-9970904, and third author by NSF
grant DMS-0070190. 1991 Mathematics Classification: 53C25, 53C12, 14E30}
\bigskip

\centerline{\vbox{\hsize = 5.85truein
\baselineskip = 12.5truept
\eightrm
\noindent {\bsc Abstract:}
On simply connected five manifolds \Se metrics coincide with Riemannian
metrics admitting real Killing spinors which are of great interest as models
of near horizon geometry for three-brane solutions in superstring theory [KW].
We expand on the recent work of Demailly and Koll\'ar [DK] and 
Johnson and Koll\'ar [JK1] who give methods for constructing K\"ahler-Einstein
metrics on log del Pezzo surfaces. By [BG1] circle V-bundles over log del
Pezzo surfaces with K\"ahler-Einstein metrics have \Se metrics on the total
space of the bundle. Here these simply connected 5-manifolds arise as links 
of isolated hypersurface singularities which by the well known work of Smale
[Sm] together with [BG3] must be diffeomorphic to 
$\scriptstyle{S^5\#l(S^2\times
S^3)}.$ More  precisely, using methods from Mori theory in algebraic geometry
we prove the existence of 14 inequivalent \Se structures on
$\scriptstyle{S^2\times S^3}$ and infinite families of such structures on
$\scriptstyle{\#l(S^2\times S^3)}$ with $\scriptstyle{2\leq l\leq7}$. We also
discuss the moduli problem for these \Se structures. 
}} 
\tenrm

\vskip .4in
\bigskip
\baselineskip = 10 truept
\centerline{\bf \intro. Introduction}  
\bigskip
Surprisingly little is known about complete Einstein metrics on
compact 5-manifolds. Until now one could list only two constructions
of such metrics. The more recent one was developed by B\"ohm [B\"o]
who obtained cohomogeneity one complete
Einstein metrics on $S^5$ and $S^2\times S^3$. 
On the other hand, the well-known
irreducible homogeneous Einstein metrics on these two spaces are the simplest
examples of what is called a \Se structure. In fact, any complex surface whose 
metric is K\"ahler-Einstein and of positive scalar curvature admits a unique
simply connected circle bundle which is canonically Sasakian-Einstein.
Since all such del Pezzo surfaces $\calz$ admitting K\"ahler-Einstein metrics
are known [Siu,Ti1-3, TY] one gets the list and the classification statement
of the first theorem of the introduction. Purely in the context of Einstein
metrics this construction is due to Kobayashi 
[Be] whereas the realization that
$\cals$ is actually a 5-manifold with real Killing spinors and, hence,
\Se space came later in the work of Friedrich and Kath [FK]. 
More generally, the importance of Killing spinors was 
already realized in 1980 in the context of
the spectrum of that Dirac operator [F]. A
classification of Riemannian manifolds admitting real Killing spinors
was finally obtained by B\"ar [B\"a] who observed that
the real Killing spinors on $M$ correspond to the parallel spinors
on $\calc(M)=(\bbr^+\times M,dt^2+t^2g)$ the metric cone on $M$. 
This simple fact allowed for an elegant 
formulation of the problem in the powerful language of holonomy groups.

The physicists' interest in five and seven-dimensional manifolds 
admitting real Killing spinors dates back to the early eighties,
where Kaluza-Klein models played a central role in
the supergravity theory. Today we witness a renewed interest in
these manifolds in the context of $p$-brane solutions in superstring theory.
These so-called $p$-branes, ``near the horizon" are modeled by
the pseudo-Riemannian geometry of the product ${\rm adS}_{p+2}\times M$, where
${\rm adS}_{p+2}$ is the $(p+2)$-dimensional anti-de-Sitter space
(a Lorentzian version of a space of constant sectional curvature) and
$(M, g)$ is a Riemannian manifold of dimension $d=D-p-2$. Here
$D$ is the dimension of the original supersymmetric theory.
In the most interesting cases of M2-branes, M5-branes, and D3-branes
$D$ equals to either 11 (M$p$-branes of M-theory) or 10 (D$p$-branes
in type IIA or type IIB string theory).
String theorists are particularly interested in
those vacua of the form ${\rm adS}_{p+2}\times M$
that preserve some residual supersymmetry. It turns out that this requirement
imposes constraints on the geometry of the Einstein manifold
$M$ which is forced to admit real Killing spinors.
Depending on the dimension $d$, the possible
geometries of $M$ are as follows: nearly K\"ahler for $d=6$, weak $G_2$
holonomy for $d=7$, Sasakian-Einstein for $d=2k+1$, and
3-Sasakian for $d=4k+3$ [AFHS, MP]. Furthermore, 
given a $p$-brane solution of the above type, the interpolation between
${\rm adS}_{p+2}\times M$ and $\bbr^{p,1}\times \calc(M)$ 
leads to a conjectured duality between supersymmetric background of the form
${\rm adS}_{p+2}\times M$ and a $(p+1)$-dimensional superconformal
field theory of $n$ coincident $p$-branes located at the conical singularity
of the $\bbr^{p,1}\times \calc(M)$ vacuum. This is a generalized version of
the Maldacena's conjecture. 

In the case of D3-branes of string theory the relevant near horizon geometry
is that of ${\rm adS}_{5}\times M$, where $M$ is a \Se 5-manifold. The
D3-brane solution interpolates between ${\rm adS}_{5}\times M$ and 
$\bbr^{3,1}\times \calc(M)$, where the cone $\calc(M)$ is a
Calabi-Yau threefold. In its original version the Maladacena conjecture
(also known as AdS/CFT duality) states that the
't Hooft large $n$ limit of $N=4$ supersymmetric Yang-Mills theory
with gauge group $SU(n)$ is dual to type IIB superstring theory
on ${\rm adS}_{5}\times S^5$ [Ma, Wi]. This conjecture was further examined by 
Klebanov and Witten [KW] for the type IIB theory on
${\rm adS}_{5}\times T^{1,1}$, where $T^{1,1}$ is the other homogeneous
\Se 5-manifold $T^{1,1}=S^2\times S^3$ and the
Calabi-Yau 3-fold $\calc(T^{1,1})$ 
is simply the quadric cone in $\bbc^4$. Using the well-known fact
that $\calc(T^{1,1})$ is a K\"ahler quotient of $\bbc^4$ (or,
equivalently, that $S^2\times S^3$ is a Sasakian-Einstein
quotient of $S^7$), a dual super Yang-Mills
theory was proposed, representing D3-branes at the  conical singularities.
In the framework of D3-branes and the AdS/CFT duality the question of
what are all the possible near horizon geometries $M$ and $\calc(M)$
might be of importance. Until this year, the only known
examples of such geometries were few and they were
exactly the circle bundles over
the del Pezzo surfaces (see the comment of Yau in [Y]). 
This has drastically changed now as we shall show in our paper. 

Recently Demailly and Koll\'ar have developed some new techniques to study
the existence of K\"ahler-Einstein metrics on compact Fano orbifolds [DK], and
applied their methods to prove the existence of K\"ahler-Einstein metrics on
three log del Pezzo surfaces. Following this, two of the present authors [BG3]
showed how one can use the results of [DK] to obtain new \Se 5-manifolds.
Later Johnson and Koll\'ar [JK1] discovered a more efficient way to treat the
algebraic equation involved, and with the aid of a computer program gave many
more examples of log del Pezzo surfaces in weighted projective 3-spaces,
including one infinite series example. All of the [JK1] examples have Fano
index equal to one.  In the current paper we extend the results of [JK1]
to the case of higher index as well as use these results together with the
results of [JK1] to construct a plethora of \Se metrics in dimension five.
This is accomplished by realizing the compact simply connected 5-manifolds as
links of isolated hypersurface singularities given by weighted homogeneous
polynomials in $\bbc^4$. The projectivization then gives the log del Pezzo
surfaces $\calz$ and proving the existence of a K\"ahler-Einstein metric
$\calz$ is then tantamount to proving the existence of an \Se metric on the
5-manifold. It turns out that from the work of Smale [Sm] one can deduce that
these 5-manifolds are all of the form $\cals_l=\#l(S^2\times S^3)$ for $1\leq
l\leq 7.$ We also show that the moduli of \Se structures on $\#l(S^2\times
S^3)$ is exceedingly rich. First it is known [FK,BG1] that the regular \Se
structures on $\#l(S^2\times S^3)$ satisfy

\noindent{\sc Theorem} [FK,BG1]: \tensl Let $\cals_l=S^5\#l(S^2\times S^3).$
\item{1)} For each $l=0,1,3,4,$ there is precisely one regular \Se structure
on $\cals_l.$
\item{2)} For each $5\leq l\leq 8$ there is a $2(l-4)$ complex parameter family of
inequivalent regular \Se structures on $\cals_l.$
\item{3)} For $l=2$ or $l\geq 9$ there are no regular \Se structures on 
$\cals_l.$
\tenrm

The first non-regular \Se structures on any 5-manifold were given in [BG3].
Here the two of the present authors showed the existence of two inequivalent
non-regular \Se structures on $S^2\times S^3,$ and one non-regular \Se
structure on $\#2(S^2\times S^3).$ The latter was the first example of any \Se
metric on $\#2(S^2\times S^3).$ In the present work we use the methods
presented in [BG3] as well as those in [JK1] to greatly expand the list of
non-regular \Se structures on simply connected 5-manifolds. Explicitly we
prove the  following:

\noindent{\sc Theorem} A: \tensl 
\item{1)} $S^2\times S^3$ admits $14$ inequivalent 
non-regular Sasakian-Einstein
structures.
\item{2)} $\#2(S^2\times S^3)$ admits two distinct 1-complex parameter
families plus $21$ inequivalent non-regular Sasakian-Einstein structures.
\item{3)} $\#3(S^2\times S^3)$ admits two distinct 2-complex parameter
families, and four distinct 1-complex parameter families, and one countably
infinite family of inequivalent non-regular Sasakian-Einstein
structures. 
\item{4)} $\#4(S^2\times S^3)$ admits two distinct 3-parameter
complex families, one 2-complex parameter family, two countably infinite
complex families, and two distinct inequivalent non-regular Sasakian-Einstein
structures. 
\item{5)} $\#5(S^2\times S^3)$ admits one 4-complex parameter
family, one 3-complex parameter family, and two countably infinite
families of inequivalent non-regular Sasakian-Einstein structures.     
\item{6)} $\#6(S^2\times S^3)$ admits one 5-complex parameter family, one
3-complex parameter family, and two countably infinite families of
inequivalent non-regular Sasakian-Einstein structures. 
\item{7)} $\#7(S^2\times S^3)$ admits a countably infinite series of
5-complex parameter families of inequivalent non-regular Sasakian-Einstein
structures. 

\noindent Furthermore, inequivalent Sasakian structures correspond to
inequivalent Riemannian metrics. \tenrm

In light of these results an outstanding question is:

\noindent{\sc Question}: \tensl Are there any \Se structures on $\#l(S^2\times
S^3)$  for $l\geq 9?$ \tenrm

It is well-known for the regular case
that the second Betti number of any smooth del Pezzo surface is less than or
equal to 9. However, in stark contrast to this one can construct log del
Pezzo surfaces with any desired second Betti number, but as shown in Theorem
\alg.5 below the methods used in this paper for proving the existence of
K\"ahler-Einstein metrics apply only to log del Pezzo surfaces with second
Betti number less than or equal to 10.  Indeed, we shall show
in a forthcoming work that our methods prove the existence of a 7-complex
parameter family of K\"ahler-Einstein metrics on a log del Pezzo surface whose
second Betti number is 10; hence, proving the existence of a 7-complex
parameter family of \Se structures on $\#9(S^2\times S^3).$ 

\bigskip
\noindent{\sc Acknowledgments}: We would like to thank Alex Buium,
for several discussions, and  J\'anos Koll\'ar for many valuable e-mail
communications as well as his interest in our work. We also thank Gang Tian for
e-mail communications, and Gueo Grantcharov for telling us about Pontecorvo's
result [Pon] which allowed us to sharpen some results on the moduli
spaces which appeared in the first version of this manuscript.   
We owe special thanks to Jennifer Johnson for letting us use the
original C program employed in [JK1]. Finally, we would like to thank  Evan
Thomas for helping us modify this program and make it more efficient in our
calculations for the higher index case.

\bigskip
\baselineskip = 10 truept
\centerline{\bf \trans. The Transverse Geometry of a Sasakian Manifold}  
\bigskip

In this section we study the transverse geometry of the Riemannian
foliation $\calf_\xi$ of a Sasakian manifold $M.$  Good references for the
transverse geometry of foliations are [Ton] and [Mol]. We first make note of
some well known properties. The foliation $\calf_\xi$ is one dimensional whose
leaves are geodesics with respect to the Sasakian metric $g,$ and this metric
is bundle-like.  

Let $(M,\xi,\eta,\Phi,g)$ be a Sasakian manifold, and consider the contact
subbundle $\cald=\ker~\eta.$ There is an orthogonal splitting of the tangent
bundle as
$$TM=\cald \oplus L_\xi, \leqno{\trans.1}$$
where $L_\xi$ is the trivial line bundle generated by the Reeb vector field
$\xi.$ The contact subbundle $\cald$ is just the normal bundle to the
characteristic foliation $\calf_\xi$ generated by $\xi.$ It is naturally
endowed with both a complex structure $J=\Phi|\cald$ and a symplectic
structure $d\eta.$ Hence, $(\cald,J,d\eta)$ gives $M$ a {\it transverse
K\"ahler} structure with K\"ahler form $d\eta$ and metric $g_\cald$ defined by
$$g_\cald(X,Y)=d\eta(X,JY) \leqno{\trans.2}$$
which is related to the Sasakian metric $g$ by
$$g=g_\cald \oplus \eta\otimes \eta. \leqno{\trans.3}$$
Recall [Ton] that a smooth p-form $\gra$ on $M$ is called {\it basic} if 
$$\xi\rfloor \gra=0, \qquad \call_\xi\gra=0, \leqno{\trans.4}$$
and we let $\grL^p_B$ denote the sheaf of germs of basic p-forms on $M,$ snd
by $\grO_B^p$ the set of global sections of $\grL^p_B$ on $M.$ The sheaf
$\grL^p_B$ is a module under the ring, $\grL^0_B,$ of germs of smooth basic
functions on $M.$ We let $C^\infty_B(M)=\grO^0_B$ denote global sections of
$\grL^0_B,$ i.e. the ring of smooth basic functions on $M.$  Since exterior
differentiation preserves basic forms we get a de Rham complex
$$\cdots\ra{2.5}\grO_B^p\fract{d}{\ra{2.5}}\grO_B^{p+1}\ra{2.5}\cdots
\leqno{\trans.5}$$  
whose cohomology $H^*_B(\calf_\xi)$ is called the
{\it basic cohomology} of $(M,\calf_\xi).$ The basic cohomology ring
$H^*_B(\calf_\xi)$ is an invariant of the foliation $\calf_\xi$ and hence, of
the Sasakian structure on $M.$ It is related to the ordinary de Rham
cohomology $H^*(M,\bbr)$ by the long exact sequence [Ton]
$$\cdots\ra{2.5}H_B^p(\calf_\xi)\ra{2.5}H^p(M,\bbr)\fract{j_p}{\ra{2.5}}
H_B^{p-1}(\calf_\xi) \fract{\grd}{\ra{2.5}} H^{p+1}_B(\calf_\xi)\ra{2.5}\cdots  
\leqno{\trans.6}$$  
where $\grd$ is the connecting homomorphism given by $\grd[\gra]_B=[d\eta\wedge
\gra]_B=[d\eta]_B\cup[\gra]_B,$ and $j_p$ is the composition of the map induced
by $\xi\rfloor$ with the well known isomorphism $H^r(M,\bbr)\approx
H^r(M,\bbr)^{S^1}$ where $H^r(M,\bbr)^{S^1}$ is the $S^1$-invariant
cohomology defined from the  $S^1$-invariant r-forms $\grO^r(M)^{S^1}.$ Here
we denote cohomology classes in $H^p_B(\calf_\xi)$ by $[\cdot]_B$ in order to
distinquish them from the ordinary cohomology classes. We also note that
$d\eta$ is basic even though $\eta$ is not. 

Next we exploit the fact that the transverse geometry is K\"ahler [ElK]. Let
$\cald_\bbc$ denote the complexification of $\cald,$ and decompose it into its
eigenspaces with respect to $J,$ that is, $\cald_\bbc= \cald^{1,0}\oplus
\cald^{0,1}.$ Similarly, we get a splitting of the complexification of the
sheaf $\grL^1_B$ of basic one forms on $M,$ namely
$$\grL^1_B\otimes \bbc = \grL^{1,0}_B\oplus \grL^{0,1}_B.$$
We let $\grL^{p,q}_B$ denote the sheaf of germs of basic forms of type
$(p,q),$ and as in the usual case there is a splitting
$$\grL^r_B\otimes \bbc = \bigoplus_{p+q=r}\grL^{p,q}_B, \leqno{\trans.7}$$ 
as well as the {\it basic Dolbeault complex}
$$0\ra{1.5}\grL^{p,0}_B \fract{\bar{\partial}}{\ra{1.5}}\grL^{p,1}_B
\fract{\bar{\partial}}{\ra{1.5}}\cdots \ra{1.5}\grL^{p,n}_B \ra{1.5} 0,
\leqno{\trans.8}$$
together with its {\it basic Dolbeault cohomology groups}
$H^{p,q}_B(\calf_\xi).$ Most of the usual results about K\"ahler geometry
carry over to transverse K\"ahler geometry [ElK].
In particular, we have the following Proposition which follows from transverse
K\"ahler geometry:

\noindent{\sc Proposition} \trans.9: \tensl Let $(M,\xi,\eta,\Phi,g)$ be a
compact Sasakian manifold of dimension $2n+1.$ Then we have
\item{(1)} $H^{n,n}_B(\calf_\xi)\approx H^{2n}_B(\calf_\xi)\approx \bbr.$
\item{(2)} The class $[d\eta]_B\in H^{1,1}_B(\calf_\xi)$ is nontrivial.
\item{(3)} $H^{p,p}_B(\calf_\xi)>0.$
\item{(4)} $H^{2p+1}_B(\calf_\xi)$ has even dimension.
\item{(5)} $H^1(M,\bbr)\approx H^1_B(\calf_\xi).$
\item{(6)} $H^r_B(\calf_\xi)= \bigoplus_{p+q=r} H^{p,q}_B(\calf_\xi).$
\item{(7)} Complex conjugation induces an anti-linear isomorphism
$H^{p,q}_B(\calf_\xi)\approx H^{q,p}_B(\calf_\xi).$
\item{(8)} If $\omega$ is a closed real $(1,1)$ form such that $\omega\in
[d\eta]_B,$ then there exists a smooth basic function $\phi$ such that $\gro
=d\eta +i\partial\bar{\partial}\phi.$
\tenrm

\noindent{\sc Proof}: (1): The exact sequence \trans.6 gives an
isomorphism $H^{2n+1}(M)\fract{j_{2n+1}}{\approx} H^{2n}_B(\calf_\xi).$ (2):
$d\eta$ is a (1,1) form by \trans.2, and since $\eta\wedge (d\eta)^n$ is a
volume form on $M,$ the isomorphism $j_{2n+1}$ implies that $d\eta$ defines a
nontrivial class in $H^{2}_B(\calf_\xi).$ (3) follows since $[d\eta]_B$ cups
to a non-trivial element in $H^{n,n}_B(\calf_\xi).$ (4): This follows from
transverse Hodge theory in the usual way [ElK]. 

To prove (5) we notice that the
beginning of the exact sequence \trans.6 is 
$$0\ra{1.7}H^1_B(\calf_\xi)\ra{1.7}H^1(M,\bbr)\fract{j_1}{\ra{1.7}}\bbr
\fract{\grd_0}{\ra{1.7}}H^2_B(\calf_\xi).$$
By (2) $\grd_0$ is injective, so
$\{0\}=\ker~\grd_0=\im~j_1;$  hence, $j_1$ is the zero map.
(6) and (7) follow from (1) and Theorem 3.4.6 of [ElK], and (8) is Proposition
3.5.1 of [ElK]. \hfill\za  

Next we define the {\it transverse Ricci tensor} $\Ric^T_g$ of $g$ 
to be the Ricci tensor of $g_\cald.$ It is related to the Ricci tensor
$\Ric_g$ of $g$ by
$$\Ric^T_g=\Ric_g|_{\cald\times \cald} + 2g|_{\cald\times \cald}. 
\leqno{\trans.10a}$$
Similarly, we define the {\it Ricci form} $\rho_g$ and {\it transverse Ricci
form} $\rho^T_g$ by  
$$\rho_g(X,Y)=\Ric_g(X,\Phi Y), \qquad \rho^T_g(X,Y) = \Ric^T_g(X,\Phi Y)
\leqno{\trans.10b}$$ 
for $X,Y$ smooth sections of $\cald.$ It is easy to check that these are
anti-symmetric of type $(1,1)$ and  are related by
$$\rho^T_g=\rho_g +2d\eta. \leqno{\trans.10c}$$

Now the contact subbundle $\cald$ is a
complex vector bundle and thus has a first Chern class $c_1(\cald)\in
H^2(M,\bbz).$ 
Consider the long exact sequence \trans.6 together with the natural map    
$H^2(M,\bbz)\ra{1.6} H^2(M,\bbr)$ whose kernel is the torsion part of
$H^2(M,\bbz).$ From (5) of Proposition \trans.9 we have
$$\matrix{&H^2(M,\bbz)&\cr
          &\decdnar{}&\cr
          &0\ra{1.2} \bbr\fract{\grd}{\ra{1.5}}
H^2_B(\calf_\xi)\fract{\gri_*}{\ra{1.5}}H^2(M,\bbr)\ra{1.6}
H^1(M,\bbr)\ra{1.6}\cdots.} \leqno{\trans.11}$$ 
As in \trans.6 the map $\grd$ is given by
$\grd(c)=c[d\eta]$ where $c\in \bbr.$  Now on a Sasakian manifold the vector bundle
$\cald^{1,0}$ is holomorphic with respect to the CR-structure, so we can
compute the free part of $c_1(\cald)=c_1(\cald^{1,0})$ from the transverse
K\"ahler geometry in the usual way. That is $c_1(\cald)$ can be represented by
a basic real closed $(1,1)$-form $\rho_B.$ The class $c_1^B=[\grr_B]\in
H^2_B(\calf_\xi)$ is independent of the transverse metric and basic connection
used to compute it, and depends only on the foliated manifold $(M,\calf_\xi)$
with its CR-structure. It is described in [ElK] and called the {\it basic
first Chern class of $\cald$} there.  Alternatively, we can think of $c_1^B$
as the negative of the first Chern class of the ``transverse canonical bundle''
$K=(\grL_B^{1,0})^n$ of $M.$ 

We now consider deformations of the Sasakian structure which fix the
basic first Chern class $c_1^B.$ We do this by considering deformations that
fix the foliation $\calf_\xi,$ in fact fix the characteristic vector field
$\xi.$ Let  $\eta_t$ be a continuous one parameter family of real1-forms
obtained by adding to $\eta$ a continuous family of basic 1-forms $\grz_t$ 
so that $\eta_t=\eta +\grz_t$  satisfies the conditions   
$$\eta_0=\eta, \qquad \grz_0=0,\qquad \eta_t\wedge (d\eta_t)^n\neq 0~~
\forall~~ t\in [0,1]. \leqno{\trans.12}$$   
This last non-degeneracy condition
implies that $\eta_t$ is a contact form on $M$ for all $t\in [0,1].$ Then by
Gray's Stability Theorem  $\eta_t$ belongs to the same underlying contact
structure as $\eta.$ Moreover, since $\grz_t$ is basic $\xi$ is the Reeb
(characteristic) vector field associated to $\eta_t$ for all $t.$ Now let us
define $$\eqalign{\Phi_t&=\Phi -\xi\otimes \grz_t\circ \Phi, \cr
           g_t&=d\eta_t\circ(\hbox{id}\otimes \Phi_t)+\eta_t\otimes \eta_t.}
\leqno{\trans.13}$$
In [BG3] it was proved that for each $t\in [0,1],$ $(\xi,\eta_t,\Phi_t,g_t)$
defines a Sasakian  structure on $M$ associated to the foliation $\calf_\xi$
and belonging to the same underlying contact structure as $\eta.$ 
These new Sasakian structures correspond to a different splitting of the
tangent bundle, namely
$$TM=L_\xi \oplus \cald_t \leqno{\trans.14}$$
which is orthogonal with respect to the metric $g_t,$ where $\cald_t=\ker~
\eta_t.$ 

Given a Sasakian structure $(\xi,\eta,\Phi,g)$ on a manifold $M,$ we 
define $\gF(\xi)$ to be the family of all Sasakian structures obtained by the
deformations above. We are now ready for

\noindent{\sc Definition} \trans.15: \tensl Two Sasakian
structures $\cals =(\xi,\eta,\Phi,g)$ and $\cals' =(\xi',\eta',\Phi',g')$ on
$M$ are said to be {\it homologous} if $\xi'=\xi$ and $[d\eta]_B=[d\eta']_B\in
H^{1,1}_B(\calf).$ In this case we also say that the Sasakian metrics $g$ and
$g'$ are {\it homologous}. \tenrm

Notice that if $L_{\xi'}=L_\xi$ we can always scale the Sasakian structure
by choosing $\xi'=\xi,$ so there is no loss in generality by choosing the
characterstic vector fields to coincide. We have

\noindent{\sc Proposition} \trans.16: \tensl Any two Sasakian structures in
$\gF(\xi)$ are homologous and the class  $c_1^B\in H^{1,1}_B(\calf_\xi)\subset
H^2_B(\calf_\xi)$  depends only on the family $\gF(\xi).$ \tenrm

\noindent{\sc Proof}: The first statement follows from (8) of Proposition
\trans.9, while the second statement follows from the fact that for all $t\in
[0,1]$ the complex vector bundles $(\cald_t,\Phi_t)$ are isomorphic. The
isomorphism between $\cald$ and $\cald_t$ is given by the map
$$\bbi-\xi\otimes \grz_t: TM\ra{1.5} TM,$$ 
and the induced map on the exterior bundle $\grL\cald$ is
$$\bbi-\grz_t\otimes \xi$$
which is the identity on basic forms. \hfill\za

Now we wish to consider the Sasakian analogue of the Calabi
problem. Its solution essentially follows from the `transverse Yau Theorem'
given by El Kacimi-Alaoui [ElK]: 

\noindent{\sc Theorem} \trans.17[ElK]: \tensl If $c_1^B$ is represented
by a real basic $(1,1)$ form $\grr^T,$ then it is the Ricci curvature form
of a unique transverse K\"ahler form $\gro^T$ in the same basic cohomology
class as $d\eta.$  \tenrm

This is given in local CR-coordinates
$(z_i,\bar{z}_i,x)$ on $M$ by solving the ``transverse Monge-Ampere equation''
$${\det(g^T_{i\barj}+\phi_{i\barj})\over \det(g^T_{i\barj})}= e^{-k\phi+F},
\quad g^T_{i\barj}+\phi_{i\barj}>0, \leqno{\trans.18}$$
for the $k=0$ case.  Here $g^T$ is the transverse metric, $\phi$ and $F$ are
real basic functions, and $\phi_{i\barj}$ are the components of
$i\partial\bar{\partial}\phi =d\grz_t$ with respect to the transverse
coordinates $(z_i,\bar{z}_{\barj}).$  Then by \trans.10c  this translates for
Sasakian geometry to:

\noindent{\sc Thereom} \trans.19: \tensl Let $(M,\xi,\eta,\Phi,g)$ be a
Sasakian manifold whose basic first Chern class $c_1^B$ is represented by
the real basic $(1,1)$ form $\grr,$ then there is a unique Sasakian structure
$(\xi,\eta_1,\Phi_1,g_1)\in\gF(\xi)$ homologous to $(\xi,\eta,\Phi,g)$ such
that $\grr_{g_1}=\grr-2d\eta_1$ is the Ricci form of $g_1,$ and
$\eta_1=\eta + \grz_1,$ with $\grz_1={1\over 2}d^c\phi.$ The metric $g_1$ and
endomorphism $\Phi_1$ are then given by \trans.12. \tenrm  

\noindent Notice that the Ricci forms are related by
$$\grr=\grr_g+{1\over 2}dd^c(F-\phi).\leqno{\trans.20}$$

Next we discuss a positivity requirement on the basic Chern
class of a compact Sasakian manifold. 

\noindent{\sc Definition} \trans.21: \tensl A Sasakian manifold $M$ is said to
be {\it positive} if its basic first Chern class $c_1^B$ can be
represented by a positive definite $(1,1)$-form. \tenrm

Notice that a positive Sasakian manifold does not necessarily have a metric
with positive Ricci curvature; however, Theorem \trans.19
does imply:

\noindent{\sc Proposition} \trans.22: \tensl A Sasakian manifold
$(M,g,\xi,\eta,\Phi)$ is positive if and only if there is a Sasakian metric
$g'$  homologous to $g$ whose Ricci curvature
satisfies the bound  $\Ric_{g'}> -2.$ In particular, a complete positive
Sasakian manifold is compact with finite fundamental group. \tenrm 

Since a \Se manifold necessarily has positive Ricci tensor, its 
Sasakian structure is necessarily positive. We are interested in sufficient
conditions on a Sasakian manifold that guarentee the existence of a \Se
structure. These conditions are algebraic geometric in nature, so we need to
impose another condition on the Sasakian structure, namely that it is
quasi-regular. Recall that the  {\it toral rank} or just {\it rank} defined
[BG2] is the dimension of the closure of the characteristic foliation and
denoted by $\hbox{rk}(M).$ For a Sasakian manifold $M^{2n+1}$ of dimension
$2n+1$ we have $1\leq \hbox{rk}(M)\leq n+1.$ The case $\hbox{rk}(M)=1$
corresponds to the quasi-regular case, and if $\hbox{rk}(M)>1$ there are
infinitely many rank one structures that are close in an appropriate sense
[BG2]. In the remainder of this paper we are essentially interested in the
rank one case only.  Now the aforementioned algebraic geometric conditions 
which imply the existence of a K\"ahler-Einstein metric are described in the
works of Nadel [Na] and Demailly and Koll\'ar [DK]. The main point is that the
obstructions for finding a solution to the Monge-Ampere equations \trans.18 
involves the non-triviality of certain {\it multiplier ideal sheaves}
associated with effective canonical $\bbq$-divisors on the space of leaves
$\calz.$ Consequently, if one can show that these multiplier ideal  sheaves
coincide with the full structure sheaf, one obtains the existence of a
K\"ahler-Einstein metric. Equivalently, and this is the approach of Johnson
and Koll\'ar [JK1], this is given in the language of Mori theory by the
appellation {\it Kawamata log terminal}, whose  precise technical definition is
given in \est.3 below. It is important to realize that these conditions are
only {\it sufficient conditions} for the existence of K\"ahler-Einstein
metrics. Indeed, they do not hold in the case of complex projective space. We
can easily reformulate the Demailly-Koll\'ar result in terms of Sasakian
geometry as

\noindent{\sc Theorem} \trans.23: \tensl Let $(\cals,g,\xi)$ be a rank 1
positive Sasakian manifold and let $\calz$ denote the space of leaves of the
characteristic foliation, and $K_\calz$ be its canonical bundle. Suppose
further that for some $\gre > 0$ and 
every effective $\bbq$-divisor $D$ on $\calz$ numerically
equivalent to $-K_\calz$ the pair $(\calz,{n+\gre\over n+1}D)$ is klt.
Then there is a \Se metric $g'$ on $M$ homologous to $g.$ \tenrm

\noindent{\sc Proof}: By Theorem 2.1 of [BG1] $\calz$ is a $\bbq$-factorial
Fano variety,and by Demailly and Koll\'ar [DK] as stated in Theorem 7 of [JK1]
$\calz$ has a K\"ahler-Einstein metric. Thus, again by (iv) of Theorem 2.1 of
[BG1] $\cals$ has a \Se metric. \hfill\za

Finally we recall [BG1] that the {\it orbifold Fano index} 
or just {\it index} of a
positive rank 1 Sasakian manifold $\cals$ is the largest positive integer $I$
such that ${1\over I}K_\calz$ is an element of $\hbox{Pic}^{orb}(\calz)$ where
$K_\calz$ is the canonical V-bundle on $\calz.$  Thus, the index is just the
generalization to the orbifold category of the usual notion of Fano index;
however, it is important to notice that the index $I$ is an invariant of the
Sasakian structure on $\cals;$ in fact, we have

\noindent{\sc Proposition} \trans.25: \tensl The orbifold Fano index $I$ is an
invariant of the family of Sasakian structures $\gF(\xi).$ \tenrm

\bigskip
\bigskip
\baselineskip = 10 truept
\centerline{\bf \est. The Existence of K\"ahler-Einstein Metrics}  
\bigskip

In this section we quickly review the basic definitions and theorems from
minimal model theory; all of this can be found in [KMM].
Throughout this section $X$ will denote a normal projective variety
defined over an algebraically closed field $k$ of characteristic zero.  By a
variety we mean an integral separated scheme of finite type over $k$. 
Since $X$ is normal,
${\rm codim}(X_{\rm sing},X) \geq 2,$ where $X_{\rm sing}$ denotes
the singular locus of $X$.  The sheaf of K\"{a}hler differentials $\Omega_X$ is
locally free of rank ${\rm dim}\,X$ on 
$X_{\rm reg} := X \backslash X_{\rm sing}$.
This gives a well-defined invertible sheaf
$\bigwedge^{{\rm dim}\,X}\Omega_{X_{\rm reg}}$ on $X_{\rm reg}$.  

\noindent{\sc Definition} \est.1: \tensl
A {\it canonical divisor} $K_X$ is any Weil divisor with
$$
{\cal O}_{X_{\rm reg}}(K_X) \simeq
 \bigwedge^{{\rm dim}\,X}\Omega_{X_{\rm reg}}.
$$
\tenrm

$K_X$ is well defined modulo linear equivalence since 
${\rm codim}(X_{\rm sing}, X) \geq 2$.  
Because singularities are essential in the study of minimal models, it is
important to note the distinction between Weil divisors and Cartier divisors.

\noindent{\sc Definition} \est.2: \tensl
A Weil divisor $D \subset X$ is called $\bbq$--{\it Cartier} if some integral
multiple $mD$ is a Cartier divisor, i.e. if ${\cal O}_X(mD)$ is invertible.
We say that $X$ is $\bbq$--{\it factorial} (or that $X$ has at worst 
$\bbq$--factorial singularities) if all Weil divisors are $\bbq$--Cartier.
\tenrm

There is a certain class of singularities called ``log--terminal'' which
play an especially prominent role in the minimal model program. 

\noindent{\sc Definition} \est.3: \tensl
Let $X$ be a normal variety and let $\Delta = \sum a_iE_i$ be an effective
$\bbq$--Weil divisor such
that $K_X + \Delta$ is $\bbq$--Cartier.  Assume $0 \leq a_i <1$ and that
the $E_i$ are all distinct.  Then the pair $(X,\Delta)$ has only
{\it log--terminal singularities} (also called klt or Kawamata
log-terminal)
if there exists a resolution of singularities
$\pi: Y \rightarrow X$ such that the union of the $\pi$--exceptional locus and
$\pi^{-1}(\cup E_i)$ is a normal crossing divisor and
$$
K_Y \equiv \pi^*(K_X + \Delta) + \sum a_iF_i
$$
with $a_i > -1$ for all $\pi$--exceptional $F_i$ (here $\equiv$ denotes
numerical equivalence).  If $a_i \geq -1$
then the pair $(X,\Delta)$ is said to have {\it log--canonical singularities}.
In the case where
$\Delta = 0$,  $X$ is said to have has log--terminal (or log--canonical)
singularities.  When
$\Delta = 0$ and $a_i > 0$ (resp. $a_i \geq 0$) $X$ is said to have
 {\it terminal singularities} (resp. {\it canonical singularities}).
\tenrm

Note that Koll\'ar and Mori, [KM] Definition 2.34,
give a different definition of klt singularities,
one which does not require the ``boundary divisor'' $\Delta$ to be effective.
However, in applications one almost always needs the effectivity hypothesis
in order to apply the machinery of the minimal model program.  All varieties
we encounter will have mild quotient singularities and hence will be
$\bbq$--factorial. 

\noindent{\sc Definition} \est.4: \tensl
Suppose $L$ is a line bundle on a variety $X$.  We write
$$
{\rm BS}(L) = \{x \in X: s(x) = 0 \ \ {\rm for\  all\  sections}
\ s \in H^0\left(X,L\right) \}
$$
for the base locus of $L$. If $\cali \subset \calo_X$ is an ideal sheaf then
we will denote by ${\rm BS}(L \otimes \cali)$ the base locus of those sections
$\sigma \in H^0(X,L)$ whose zero locus contains the zero locus of $\cali$. 
\tenrm
\medskip

We will be interested in the case where $X$ is a projective surface with 
mild quotient singularities.  We will assume that
$-K_X$ is an ample ${\bbq}$--Cartier
divisor and we will be interested in showing that for some $\epsilon > 0$
and {\it every} effective
${\bbq}$--divisor $D \in |-K_X|$, the pair $(X,{2+\epsilon\over3}D)$ is 
klt.
In particular this is certainly true if the pair
$(X,D)$ is log--canonical for all choices of $D$.  
Checking whether or not $(X,D)$
is log--canonical will be intricate at singular points of $X$ and we begin
by reducing this condition to a much simpler one on a smooth finite cover
of $X$.

\medskip

We begin by noting 
that by [KM] Proposition
4.18, the variety $X$ itself is klt at a singular point 
$x$ if and only if $X$ has a quotient singularity at $x$.  
In order to
check whether or not the pair $(X,D)$ is klt at $x$ we need to consider
a local cover
$$
\pi: \left({\bbc}^2,0\right) \rightarrow (X,x),
$$
where the map $\pi$ is given by taking the quotient via a finite
group action.  We let $U \subset {\bbc}^2$ 
denote an open subset on which $\pi$
is defined and let
$$
p: S \rightarrow U
$$
be the blow--up of $U$ at $0$ with exceptional divisor $E$.  
We will show that if
$(S,p^{-1}_\ast (\pi^\ast D)+E)$ 
is log--canonical near $E$ then
the pair $(X,D)$ is log--canonical
at $x$.
Shokurov's inversion
of adjunction ([KM] Theorem 5.50) states that the pair 
$(S,p_\ast^{-1} (\pi^\ast D)+E)$ is log--canonical near $E$ if 
$(E,p_\ast^{-1}\pi^\ast(D)|E)$ is log--canonical.
But $E \simeq {\bbp}^1$ and so
this is an easy condition to check, satisfied provided  the round--up of
$p_\ast^{-1}\pi^\ast(D)|E$ is a sum of distinct points, each with multiplicity 
one. 
\medskip 

By [KM] Lemma 2.30 and Proposition 5.20.4 together with the inversion of
adjunction, one derives

\noindent{\sc Lemma} \est.5: \tensl
Suppose that ${\rm mult}_0(\pi^\ast D) \leq 2$ 
and that $(E,p_\ast^{-1}\pi^\ast D
+ E)$ is log--canonical.  Then $(X,D)$ is log--canonical at $x$. 
\tenrm

\medskip
We sketch here the details of Lemma 2.5.
Suppose that $(S, p_\ast^{-1}(\pi^\ast D)+E)$ is log--canonical 
near $E$ and consider the following commutative diagram:

$$\matrix{&&Y&&\cr
&\fract{q}{\swarrow}&&\fract{q'}{\searrow}&\cr
S&&\decdnar{\mu}&&X'.\cr
&\fract{\searrow}{p}&&\fract{\swarrow}{p'}&\cr
&&U&&\cr}$$

Here we choose $\mu$ so that $\mu^{-1}(\pi^\ast D) \cup {\rm Exc}(\mu)$ 
is a divisor with normal crossings.
Suppose we write
$$K_Y\equiv
q^\ast(K_S + p_\ast^{-1}(\pi^\ast D)+E) + \sum_{i=1}^s a_iE_i .
\leqno{\est.6}$$

Then $(S,p_\ast^{-1}(\pi^\ast D)+E)$ is log--canonical near $E$ provided
$a_i \geq -1$ for all  $E_i$ contracted by $q$.
We also have, using $q^\ast K_S - \mu^\ast K_U = q^\ast E$ ,

$$\leqalignno{
  K_Y&\equiv \mu^\ast(K_U + \pi^\ast(D)) + \sum_{j=1}^s b_jE_j\cr
      &\equiv q^\ast(K_S -E + p^\ast \pi^\ast D) + \sum_{j=1}^s b_jE_j\cr
      &\equiv q^\ast(K_S+ p_\ast^{-1}(\pi^\ast D)
       + E) - (2-{\rm mult}_0(\pi^\ast D))
          q^\ast E + \sum_{j=1}^s b_jE_j. &\est.7\cr} 
$$

Combining \est.6 and \est.7 shows that 
$b_i \geq a_i$ 
for all $i$, provided that ${\rm mult}_0(\pi^\ast D) \leq 2$.   
Thus, whenever ${\rm mult}_0(\pi^\ast D) \leq 2$,  $(U,\pi^\ast D)$ 
is log--canonical at $0$ provided
$(S,p_\ast^{-1} (\pi^\ast D)+E)$ is log--canonical near $E$.

\medskip
Next we need to relate the pair $(X,D)$ at $x$ 
to the pair $(U,\pi^\ast D)$ at $0$.
For this, we consider the commutative diagram

$$\matrix{&&W&&\cr
&\fract{g}{\swarrow}&&\fract{h}{\searrow}&\cr
U&&\phantom{\decdnar{\mu}}&&\tilde{X}.\cr
&\fract{\searrow}{\pi}&&\fract{\swarrow}{f}&\cr
&&X&&\cr}$$
\noindent
Write
$$\eqalign{
K_W =& g^\ast K_U + \sum_{j=1}^s b_jG_j ,\cr
K_{\tilde{X}}=& f^\ast K_X + \sum_{i=1}^r a_iF_i.\cr}$$
\noindent
Using the logarithmic ramification formula ([Ka] Lemma 1.6) we find
$$\sum_{j=1}^s(1+b_j)G_j = h^\ast \left(\sum_{i=1}^r (1+a_i)F_i\right)
+ R,\leqno{\est.8}$$
where all divisorial components of $R$ are contracted by $h$.
We now remove the boundary divisor from both sides of \est.8 giving
$$\sum_{j=1}^s(1+b_j)G_j -g^\ast \pi^\ast D = 
h^\ast \left(\sum_{i=1}^r (1+a_i)F_i - f^\ast D\right)
+R.\leqno{\est.9}$$
We rewrite \est.9 as
$$\sum_{j=1}^{s^\prime} \beta_jG_j^\prime =
h^\ast \left(\sum_{i=1}^{r^\prime} \alpha_i F_i^\prime \right)
+R.\leqno{\est.10}$$
Suppose now that $F$ is one of the $F_i^\prime$ in \est.10 which is 
$f$--exceptional and let $\alpha _F$ denote the corresponding discrepancy.  
Let $G$ denote a divisor on $W$ such that $h: G \rightarrow F$ 
is finite, with discrepancy $\beta_G$.  
Then \est.10 says that $\alpha_F \geq -1$ provided
$\beta_G \geq -1$.  
But $(U,\pi^\ast D)$ is log--canonical near $E$ if and only if $\beta_j \geq 
-1$ for all $G_j$ contracted by $g$. Thus $\alpha_F \geq -1$ for all $F$ 
contracted by $f$ and this shows that $(X,D)$ is log--canonical at $x$.

\medskip

We now give some general criteria for establishing
whether or not $(X,D)$ is log--canonical at a point $x$.

\noindent{\sc Lemma} \est.11: \tensl
Suppose $X$ is a normal projective surface and $D$ an effective
$\bbq$--divisor on $X$.  Suppose moreover that $x \in X_{\rm reg}$.  
Then the pair $(X,D)$ is klt at $x$ provided $A \cdot D < a$ for some divisor
$A$ such that the linear series $|A \otimes m_x^a|$ 
has only isolated base points.  
\tenrm

\noindent
{\sc Proof}:
Suppose $x \in {\rm supp}(D)$.  Choose a representative $E \in 
|A \otimes m_x^a|$ such that
$E \cap {\rm supp}(D)$ is proper.  Then we have
$$
i_x(D,E) \leq D \cdot E < a.
$$
On the other hand
$$
i_x(D,E) \geq {\rm mult}_x(D) {\rm mult}_x(E) \geq a({\rm mult}_x(D)).
$$
Thus we have
$$
{\rm mult}_x(D) < 1
$$
and Lemma \est.11 follows from [KM] Theorem 4.5 (1). \hfil\za  

Note that Lemma \est.11  throws away a lot of information because the
estimate
$$
D \cdot E \geq i_x(D,E)
$$
can be very crude; 
to give a more precise estimate, however, requires information
about the specific geometry of $X$.  
We will also need a slight refinement of Lemma \est.11 to deal with
the singular points of $X$.  

\noindent{\sc Lemma} \est.12: \tensl
Suppose $X$ is a normal projective surface with a quotient 
singularity of index $m$ at $x$.  Suppose $D$ is an effective
 $\bbq$--divisor and $A$ a 
divisor on $X$ such that $|A \otimes m_x^a|$ has only isolated base points for
some positive integer $a$.  Then $(X,D)$ is klt at $x$ provided $A \cdot D
< {a\over m}$.  
\tenrm

\noindent{\sc Proof}:  Choose a general divisor $E \in |A \otimes
 m_x^a|$
so that $E \cap D$ is proper.  Choose a local resolution
 $\pi: Y \rightarrow X$ of $X$ at $x$ with $\pi^{-1}(x) = y$.  Then 
$$
{\rm mult}_y(\pi^\ast D) \leq {m D \cdot A\over a}.
$$
Thus, by hypothesis ${\rm mult}_y(\pi^\ast D) < 1$ so that $(Y,\pi^\ast D)$ is 
klt at $y$.  By Lemma \est.11 this implies that the pair $(X,D)$ is
klt at $x$.\hfill\za 

\settabs 10\columns
\bigskip
\bigskip
\baselineskip = 10 truept
\centerline{\bf \link. The Sasakian Geometry of Links of Weighted
Homogeneous Polynomials}   \bigskip

In this section we discuss the Sasakian geometry of links of isolated
hypersurface singularities defined by weighted homogeneous polynomials.
Consider the affine space $\bbc^{n+1}$ together with a weighted
$\bbc^*$-action given by $(z_0,\ldots,z_n)\mapsto
(\grl^{w_0}z_0,\ldots,\grl^{w_n}z_n),$ where the {\it weights} $w_j$ are
positive integers. It is convenient to view the weights as the components of a
vector $\bfw\in (\bbz^+)^{n+1},$ and we shall assume that they are ordered
$w_0\leq w_1\leq \cdots\leq w_n$ and that $\gcd(w_0,\ldots,w_n)=1.$ Let $f$
be a quasi-homogeneous polynomial, that is $f\in \bbc[z_0,\ldots,z_n]$ and
satisfies $$f(\grl^{w_0}z_0,\ldots,\grl^{w_n}z_n)=\grl^df(z_0,\ldots,z_n),
\leqno{\link.1}$$
where $d\in \bbz^+$ is the degree of $f.$ We are interested in the {\it
weighted affine cone} $C_f$ defined by
the equation $f(z_0,\ldots,z_n)=0.$ We shall assume that the origin in
$\bbc^{n+1}$ is an isolated singularity, in fact the only singularity, of
$f.$ Then the link $L_f$ defined by 
$$L_f= C_f\cap S^{2n+1}, \leqno{\link.2}$$
where 
$$S^{2n+1}=\{(z_0,\ldots,z_n)\in \bbc^{n+1}|\sum_{j=0}^n|z_j|^2=1\}$$
is the unit sphere in $\bbc^{n+1},$ is a smooth manifold of dimension $2n-1.$ 
Furthermore, it is well-known [Mil] that the link $L_f$ is $(n-2)$-connected.

On $S^{2n+1}$ there is a well-known [YK] ``weighted'' Sasakian structure  
$(\xi_\bfw,\eta_\bfw,\Phi_\bfw,g_\bfw)$ which in the standard coordinates
$\{z_j=x_j+iy_j\}_{j=0}^n$ on $\bbc^{n+1}=\bbr^{2n+2}$ is determined by
$$\eta_\bfw = {\sum_{i=0}^n(x_idy_i-y_idx_i)\over\sum_{i=0}^n
w_i(x_i^2+y_i^2)}, \qquad \xi_\bfw
=\sum_{i=0}^nw_i(x_i\partial_{y_i}-y_i\partial_{x_i}),$$
and the standard Sasakian structure $(\xi,\eta,\Phi,g)$ on $S^{2n+1}.$
The embedding $L_f\hookrightarrow S^{2n+1}$ induces a Sasakian structure on
$L_f$ [BG3]. 

Given a sequence $\bfw =(w_0,\ldots,w_n)$ of ordered positive integers one can
form the graded polynomial ring $S(\bfw)=\bbc[z_0,\ldots,z_n]$, where $z_i$ has
grading or {\it weight} $w_i.$ The weighted projective space [Dol, Fle]
$\bbp(\bfw)=\bbp(w_0,\ldots,w_n)$ is defined to be the scheme
$\hbox{Proj}(S(\bfw)).$  It is the quotient space 
$(\bbc^{n+1}-\{0\})/\bbc^*(\bfw)$, where $\bbc^*(\bfw)$ is the weighted action
defined in \link.1, or equivalently, $\bbp(\bfw)$ is the quotient of
the weighted Sasakian sphere
$S_\bfw^{2n+1}=(S^{2n+1},\xi_\bfw,\eta_\bfw,\Phi_\bfw,g_\bfw)$ by the
weighted circle action $S^1(\bfw)$ generated by $\xi_\bfw.$ As such
$\bbp(\bfw)$ is also a compact complex orbifold with an induced K\"ahler
structure. We have from [BG3]

\noindent{\sc Theorem} \link.3: \tensl The quadruple
$(\xi_\bfw,\eta_\bfw,\Phi_\bfw,g_\bfw)$ gives $L_f$ a quasi-regular Sasakian
structure such that there is a commutative diagram
$$\matrix{L_f &\ra{2.5}& S^{2n+1}_\bfw&\cr
  \decdnar{\pi}&&\decdnar{} &\cr
   \calz_f &\ra{2.5} &\bbp(\bfw),&\cr}$$
where the horizontal arrows are Sasakian and K\"ahlerian embeddings,
respectively, and the vertical arrows are principal $S^1$ V-bundles and
orbifold Riemannian submersions.   Moreover, $L_f$ is the total space
of the principal $S^1$ V-bundle over the orbifold $\calz_f$ whose first Chern
class  in $H^2_{orb}(\calz_f,\bbz)$ is $c_1(\calz_f)/I,$ where $I$ is the
index. \tenrm 

We should also mention that 
$c_1(\calz_f)$ pulls back to the basic first Chern
class $c_1^B\in H^2_B(\calf_{\xi_\bfw})$ and $\eta_\bfw$ is the connection in
this V-bundle whose curvature is $d\eta ={2ni\over I}\pi^*\gro_\bfw,$ where
$\gro_\bfw$ is the K\"ahler form on $\calz_f.$ 

Now conditions on the weights that guarantee that the hypersurface $C_f\subset
\bbc^{n+1}$ have only an isolated singularity at the origin are well-known
[Fle,JK1]. These conditions become more complicated as the dimension increases
[Fle,JK2]; however, in this paper we will only be interested in the $n=3$ case
of hypersurfaces in a weighted complex projective 3-space.  These conditions,
known as {\it quasi-smoothness} conditions gaurentee that $\calz_f$ is smooth
in the orbifold sense, that is, at a vertex  $P_i\in\bbp({\bf w})$ the preimage
of $\calz_f$ in the orbifold chart of  $\bbp({\bf w})$ is smooth.  It is
easy to see that one can formulate all these conditions as follows [Fle,JK1]: 
\bigskip
\noindent{\sc Quasi-Smoothness Conditions} \link.4:
\bigskip
\+I.&For each $i=0,\cdots,3$ there is a $j$ and a monomial
$z_i^{m_i}z_j\in \calo(d).$\cr
\+&Here $j=i$ is possible.\cr
\medskip
\+II.&If $\gcd(w_i,w_j)>1$ then there is a 
monomial $z_i^{b_i}z_j^{b_j}\in
\calo(d).$\cr
\medskip
\+III.&For every $i,j$ either 
there is a monomial $z_i^{b_i}z_j^{b_j}\in
\calo(d),$\cr
\+&or there are monomials $z_i^{c_i}z_j^{c_j}z_k$ and
$z_i^{d_i}z_j^{d_j}z_l\in \calo(d)$ with $\{k,l\}\neq \{i,j\}.$\cr

In condition I the $i=j$ case corresponds to the
case when $\calz_f$ does not pass through the point $P_i$. The second
condition is equivalent to $\calz_f$ not containing any of the singular
lines in $\bbp({\bf w})$. If $\calz_f$ contains a
coordinate axis (say $z_j=z_j=0$) then  the condition III
forces $\calz_f$ to be smooth
along it, except possibly at the vertices.

There is another condition apart from quasi-smoothness that assures us
that the adjunction theory behaves correctly, and that $\bbp({\bf w})$ 
does not have any orbifold singularities of codimension 1. It is [Dol,Fle] 
\medskip
\noindent{\sc Well-formedness Condition} \link.5
\bigskip
\+IV.& For any triple $i,j,k\neq,$ we have $\gcd(w_i,w_j,w_k)=1.$ \cr
\medskip

Condition IV
guarantees that the canonical V-bundle $K_\calz$ is
determined in terms of the degree and index by
$$K_\calz \simeq \calo(-I)=\calo(d-|\bfw|),\leqno{\link.6}$$
where $|\bfw|=\sum_iw_i.$ 

Finally we end this section by giving a corollary of Lemma \est.12 which, by
Theorem \trans.23, gives sufficient conditions 
to ensure that $\calz_f$ admits a
K\"ahler-Einstein metric.

\noindent{\sc Corollary} \link.7: \tensl Let  $\calz_f \subset
{\bbp}(w_0,w_1,w_2,w_3)$ be a hypersurface  of degree $d$ in weighted
projective space with $\bfw=(w_0,w_1,w_2,w_3)$ 
well-formed.  Let $\Delta \in  |-\alpha K_{\calz_f}|$ be
an effective representative for some rational $\alpha >0$. 
Writing $K_{\calz_f} = \calo_{\calz_f}(-I)$,  
suppose $x\in \calz_f$ is a point
at which ${\bbp(\bfw)}$ has a singularity of order $\ell_x.$  
If ${\rm BS}(\calo_{\calz_f}(w_0) \otimes
\cali_x) $ does not contain any component of $\Delta$
then $(\calz_f,\Delta)$ is klt
at $x$ provided
$$
\alpha \ell_x d I < w_1w_2w_3.
$$
If ${\rm BS}(\calo_{\calz_f}(w_0w_1) \otimes \cali_x)$ does not contain
$\Delta$ then  $(\calz_f,\Delta)$ is
klt at $x$ provided
$$
\alpha \ell_x d I  <  w_0w_2w_3 .
$$ 
In general, unless $\calz_f$ is the hypersurface $\{z_0 = 0\}$,
 one always has that $(\calz_f,\Delta)$ is klt
at $x$ if
$$
\alpha \ell_x d I < w_0w_1w_3.
$$ 
\tenrm

\noindent{\sc Proof}:
This corollary follows immediately from Lemma \est.12 by looking
at the linear systems $|\calo_{\calz_f}(w_0)|, |\calo_X(w_0w_1)|$, and 
$|\calo_{\calz_f}(w_0w_1w_2)|$ respectively.  For example, taking $A = 
|\calo_{\calz_f}(w_0w_1)|$
in Lemma \est.12, we find that we can take $a = w_0$ since one of
the sections $z_0^{w_1}$, $z_1^{w_0}$ will not vanish along $\Delta$ 
by hypothesis and $w_0 \leq w_1$.  
Then Lemma \est.12 states that $(\calz_f,\alpha \Delta)$ is klt at $x$ provided
$A \cdot \alpha \Delta < {w_0\over \ell_X}$.
But 
$$
A \cdot \alpha\Delta = 
{w_0w_1d \alpha I\over w_0w_1w_2w_3} = {d  
\alpha I\over
w_2w_3}.
$$
The middle formula follows immediately and the others are obtained similarly.
\hfil\za

\def\ss{\scriptscriptstyle}
\bigskip
\bigskip
\baselineskip = 10 truept
\centerline{\bf \alg. The Algebraic Equations and Their Solutions}  
\bigskip

We want to repeat the analysis of [JK1] for the case of an arbitrary
index $I$. Following their approach we consider the linear system
$$m_iw_i+w_{j(i)}=d=w_0+w_1+w_2+w_3-I\leqno{\alg.1}$$
where $m_i$ is a positive integer and both $i$ and $j(i)$ are integers ranging
over $0,1,2,3.$ Note we can have $j(i)=i.$ This is an obvious translation of
the quasi-smoothness condition $\link.4.I$ 

\noindent{\sc Lemma} \alg.2: \tensl Assuming conditions \link.4.I and
\link.4.IV, the following bounds hold 
\item{(1)} $1\leq m_3\leq 2.$
\item{(2)} Either $2\leq m_2\leq 4,$ or $2w_0\leq I.$ 
\item{(3)} Either $2\leq m_1\leq 10,$ or $2I=w_0+w_1,$ or $2I\geq 3w_0,$ or
condition \link.4.II is violated. \tenrm

\noindent{\sc Proof}: Since the weights are ordered $w_0\leq w_1\leq w_2\leq
w_3,$ (1) follows immediately from the linear system with $i=3.$ To prove (2)
we consider $i=1$ in the system \alg.1. There are two cases: $w_{j(2)}=w_3$ and
$w_{j(2)}\neq w_3.$ If $w_{j(2)}=w_3$ we immediately have 
$$m_2w_2=w_0+w_1+w_2-I\leq 3w_2-I$$
which implies $m_2\leq 2.$ 

Now consider $w_{j(2)}\neq w_3,$ then there are two subcases depending on
whether $m_3=2$ or $1.$ If $m_3=2$ our system gives
$$w_0+w_1+w_2=w_3+w_{j(3)}+I,$$
and here there are two subcases giving: 
\item{1.} If $w_{j(3)}=w_3,$ then $2w_3<3w_2.$
\item{2.} If $w_{j(3)}\neq w_3,$ then $w_3<2w_2,$ 

\noindent and in both cases we have $w_3<2w_2.$ Thus since $w_{j(2)}\neq w_3,$ we
get  
$$m_2w_2\leq 2w_2+w_3-I< 4w_2-I$$
implying $m_2\leq 3.$

Finally we consider the case $m_3=1.$ Here we obtain
$$w_{j(3)}+I=w_0+w_1+w_2.\leqno{\alg.3}$$
If $w_3=w_{j(3)}=w_0+w_1+w_2-I,$ then $d=2w_3,$
so $m_2w_2+w_{j(2)}= 2w_0+2w_1+2w_2-2I.$ This gives
$$m_2w_2\leq 5w_2-2I,$$
implying that $m_2\leq 4,$ whereas if $w_{j(3)}\neq w_3,$ then we have
$I=w_i+w_j$ for some distinct $i,j=0,1,2.$  
We have thus shown that either $1\leq m_2\leq 4$ or $I=w_i+w_j$ for  some
distinct $i,j=0,1,2.$  It is now easy to see that if $m_2=1$ then we must
have $2w_0\leq I.$ Thus, in all cases (2) of the lemma holds.   

Finally we consider (3).  We proceed along the lines of (2).  If possible we
solve the system \alg.1 for $i=2,3$ for $w_3$ and $w_2;$  however, there are
exceptional cases, that is, values of $m_3$ and $m_2$ such $w_3$ and/or $w_2 $
are free. We shall show that in these cases when $m_1$ is unbounded that either
$2I\geq 3w_0,$ or $2I=w_0+w_1,$ or condition \link.4.II is violated. Moreover,
when we can solve for $w_2$ and $w_3$ and substitute into \alg.1 with $i=1$ we
obtain an equation for $m_1$ in terms of $m_2,m_3,w_1,w_2,w_3.$ This gives a
bound for  $m_1$ which is certainly maximal for $I=1$ (in fact for $I=0$) and
one can see that $m_1\leq 10$ as in [JK].  To see the lower bound on $m_1$ we
show that if $m_1=1$ then $I\geq 2w_0.$ There are two cases to consider,
$w_{j(1)}=w_1$ and $w_{j(1)}\neq w_1,$ and one easily sees that both cases
give the desired estimate.

We now consider the exceptional cases. For the $i=3$  equation of the system
\alg.1 this occurs for $m_3=1,$ and one easily sees that there are 3
exceptional cases all which imply the estimate $I\geq 2w_0.$ Now assuming that
$w_3$ is determined in terms of the remaining weights and index, we consider
the $i=2$ equation in \alg.1. There are 4 subcases.  Since $m_2>1$ the case
$w_{j(2)}=w_3$ gives no exceptional solution.  Now if $w_{j(2)}=w_2$ we
have $m_2w_2= w_0+w_1+w_3-I$ and we must plug in the possibilities for $w_3$
determined by the $i=3$ equation. But since all of these possibilities are of
the form $w_3= aw_2+b$ with $a=0$, $a= {1\over 2}$ or $a=1,$ 
we see that there are no
exceptional solutions in this case either.  

Next we put $w_{j(2)}=w_1,$ and consider the 5 possibilities for $w_3.$ The
first is $w_3=w_0+w_1+w_2-I,$ and this leads to the exceptional case with
$m_2=2$ and $2I=2w_0+w_1\geq 3w_0.$ The next case with $2w_3= w_0+w_1+w_2-I$
is not exceptional, and the third with $w_3=w_0+w_1-I$ again gives $2I\geq
3w_0.$  The fourth case $w_3=w_0+w_2-I$, however, gives the exceptional case
$m_2=2$ and $I=w_0,$ and this requires a bit more analysis. Here we have that
$w_3=w_2$ and so the weights are $(I,w_1,w_2,w_2)$ with degree $d=2w_2+w_1.$
But this fails the quasi-smoothness condition (2) of section \link~ the $w_2$
must divide $d,$ and this can occur only if $w_2=w_1$ which violates the
well-formedness condition.  Finally the fifth subcase with $w_3=w_1+w_2-I$
leads to $m_2=2$ and $2I=w_0+w_1.$

Finally we have the case $w_{j(2)}=w_0.$ The first possibility for $w_3$ gives
an exceptional case with $2I\geq 3w_0,$ while the next two do not lead to any
exceptional cases. Subcase 4, however, with $w_3=w_0+w_2-I$ give $m_2=2$ and
$2I=w_0+w_1.$ This implies the inequalities $w_0\leq I\leq w_1$ and the
equation $w_3=w_2-w_1+I$ which violate the order $w_2\leq w_3$ unless
$w_3=w_2$ and $I=w_1=w_0.$ But in this case the degree is $d=2w_2+I,$ and
condition \link.4.II then implies that $w_2=w_3=I$ which 
contradicts the well-formedness condition \link.4.IV. The last possibility for
$w_3$ is $w_3=w_1+w_2-I$ and this leads to the exceptional case $m_2=2$ and
$I=w_1.$ But then we have $w_3=w_2$ and $d=2w_2+w_0$ which violates condition
\link.4.II above again since $w_2$ cannot divide $w_0.$ \hfill\za

\noindent{\sc Remarks} \alg.4: There are solutions with arbitrary $m_1$ 
that are quasi-smooth, but our lemma implies that these must satisfy $2I\geq
3w_0$ or $2I=w_0+w_1.$ As we shall see in the next section,
in both cases the sufficient condition for the
existence of a K\"ahler-Einstein metric fails.  Examples of such a series are:
$(1,1,k,k)$ of degree $2k$ and index 2, $(1,2,2k+1,2k+1)$ of degree
$2(2k+1)$ and index 3, and $(2,2,2k+1,2k+1)$ of degree $2(2k+1)$ and index 4,
as well as the general series $(I-n,I+n,w,w+n)$ of degree $2w+n+I$ and index
$I.$ There are also double series such as $(1,1,m,m+k)$ of degree $2m+k$ and
index 2. 

We have modified the computer program of [JK1] to solve the system \alg.1 for 
any index $I$ as well as to discard some solutions that are not
quasi-smooth.  In light of Remarks \alg.4,
we shall now consider all solutions
for which $1\leq I\leq10$, $2I< 3w_0,$ and $2I\not=w_0+w_1$. The last
column of the tables indicate whether the klt condition (which implies that
$\calz_\bfw$ admits a K\"ahler-Einstein metric) holds (Y) or is unknown (?).
The computer search gives the following complete list:

\noindent{\sc Theorem} \alg.5: \tensl Let $\calz_{\bf w}$ be a well-formed
quasi-smooth log del Pezzo surface of index $I\leq 10$ and degree $d$ embedded
in the weighted projective space $\bbp(\bfw)= \bbp(w_0,w_1,w_2,w_3).$ If
$2I\geq3w_0$   or $2I= w_0+w_1$ then for every $\epsilon > 0$ there exists
a divisor $D \in |-K_{\calz_\bfw}|$ such that
$(\calz_\bfw,{2 + \epsilon \over 3}D)$ is not klt. If
$2I<3w_0$  and $2I\neq w_0+w_1$ then
$\calz_\bfw$ must belong to one of the
following cases:

\item{(1)} If $I=1$ [JK1] then $\calz_{\bf w}$  is either one of the 
series solution below
\bigskip
\centerline{
\vbox{\tabskip=0pt \offinterlineskip
\def\tablerule{\noalign{\hrule}}
\halign to300pt {\strut#& \vrule#\tabskip=1em plus2em&
     \hfil#& \vrule#& \hfil#& \vrule#& \hfil#& \vrule#& 
     \hfil#& \vrule#\tabskip=0pt\cr\tablerule
&&\omit\hidewidth ${\bf w}$\hidewidth&&
\omit\hidewidth $d$\hidewidth&&
\omit\hidewidth $b_2$\hidewidth&&
\omit\hidewidth K-E\hidewidth&\cr\tablerule
&&$(2,2k+1,2k+1,4k+1)$&&$8k+4$&&8&&Y&\cr\tablerule}}
}
or it is one of the 19 sporadic solutions listed in Table 1.

\item{(2)} If $I=2$ then $\calz_{\bf w}$ 
is one of the 6 infinite series solutions with ${\bf w}$ and $d$ given by: 
\bigskip
\centerline{
\vbox{\tabskip=0pt \offinterlineskip
\def\tablerule{\noalign{\hrule}}
\halign to300pt {\strut#& \vrule#\tabskip=1em plus2em&
     \hfil#& \vrule#& \hfil#& \vrule#&  \hfil#& \vrule#&
     \hfil#& \vrule#\tabskip=0pt\cr\tablerule
&&\omit\hidewidth ${\bf w}$\hidewidth&&
\omit\hidewidth $d$\hidewidth&&
\omit\hidewidth $b_2$\hidewidth&&
\omit\hidewidth K-E\hidewidth&\cr\tablerule
&&$(4,2k+1,2k+1,4k)$&&$8k+4$&&7&&Y&\cr
\tablerule
&&$(3,3k+1,6k+1,9k+3)$&&$18k+6$&&6&&Y&\cr\tablerule
&&$(3,3k+1,6k+1,9k)$&&$18k+3$&&5&&?&\cr\tablerule
&&$(3,3k,3k+1,3k+1)$&&$9k+3$&&7&&Y&\cr\tablerule
&&$(3,3k+1,3k+2,3k+2)$&&$9k+6$&&5&&Y&\cr\tablerule
&&$(4,2k+1,4k+2,6k+1)$&&$12k+6$&&6&&?&\cr\tablerule}}
}
or it is one of the 25 sporadic solutions listed in Table 1.

\bigskip
\item{(3)} If $I=3$ then $\calz_{\bf w}$ is one of the 7 sporadic solutions
listed in Table 1.

\item{(4)} If $I=4$ then $\calz_{\bf w}$
is one of the 3 infinite series with ${\bf w}$ and $d$ given by:
\bigskip
\centerline{
\vbox{\tabskip=0pt \offinterlineskip
\def\tablerule{\noalign{\hrule}}
\halign to300pt {\strut#& \vrule#\tabskip=1em plus2em&
     \hfil#& \vrule#& \hfil#& \vrule#&  \hfil#& \vrule#&
     \hfil#& \vrule#\tabskip=0pt\cr\tablerule
&&\omit\hidewidth ${\bf w}$\hidewidth&&
\omit\hidewidth $d$\hidewidth&&
\omit\hidewidth $b_2$\hidewidth&&
\omit\hidewidth K-E\hidewidth&\cr\tablerule
&&$(6,6k+3,6k+5,6k+5)$&&$18k+15$&&5&&Y&\cr
\tablerule
&&$(6,6k+5,12k+8,18k+9)$&&$36k+24$&&3&&?&\cr\tablerule
&&$(6,6k+5,12k+8,18k+15)$&&$36k+30$&&4&&Y&\cr\tablerule}}
}
\bigskip
or it is one of the 10 sporadic solutions of Table 1.

\item{(5)} If $I=5$ then $\calz_{\bf w}$ is one of the 3 sporadic solutions
listed in Table 2.

\item{(4)} If $I=6$ then $\calz_{\bf w}$
is one of the 2 infinite series given by:
\bigskip
\centerline{
\vbox{\tabskip=0pt \offinterlineskip
\def\tablerule{\noalign{\hrule}}
\halign to300pt {\strut#& \vrule#\tabskip=1em plus2em&
     \hfil#& \vrule#& \hfil#& \vrule#&  \hfil#& \vrule#&
     \hfil#& \vrule#\tabskip=0pt\cr\tablerule
&&\omit\hidewidth ${\bf w}$\hidewidth&&
\omit\hidewidth $d$\hidewidth&&
\omit\hidewidth $b_2$\hidewidth&&
\omit\hidewidth K-E\hidewidth&\cr\tablerule
&&$(8,4k+1,4k+3,4k+5)$&&$12k+11$&&3&&?&\cr
\tablerule
&&$(9,3k+2,3k+5,6k+1)$&&$12k+11$&&3&&?&\cr\tablerule}}
}
\bigskip
or it is one of the 3 sporadic solutions of Table 1.

\item{(5)} If $I=7,8,9,10$ then $\calz_{\bf w}$ is one of the 6
sporadic solutions
listed in Table 1.
\tenrm
\bigskip
\vfill\eject

\def\ss{\scriptscriptstyle}
\centerline{
\vbox{\tabskip=0pt \offinterlineskip
\def\tablerule{\noalign{\hrule}}
\halign to380pt {\strut#& \vrule#\tabskip=1em plus2em&
     \hfil#& \vrule#& \hfil#& \vrule#&  \hfil#& \vrule#&
     \hfil#& \vrule#&
     \hfil#& \vrule#&
     \hfil#& \vrule#\tabskip=0pt\cr\tablerule
\omit&height2pt&\multispan{11}&\cr
&&\multispan{11}\hfil {\bf Table 1.} 
Sporadic Examples of $\calz_\bfw$ of Index $1\leq I\leq10$\hfil&\cr\tablerule
&&\omit\hidewidth Index\hidewidth&&
\omit\hidewidth ${\bf w}$\hidewidth&&
\omit\hidewidth Monomials of $f_{\bf w}$\hidewidth&&
\omit\hidewidth $d$\hidewidth&&
\omit\hidewidth $b_2$\hidewidth&&
\omit\hidewidth K-E\hidewidth&\cr\tablerule
&&1&&(1,2,3,5)&&$z_0^{10},z_1^5,z_2^3z_1,z_3^2,\ldots
(17)^*$&&10&&9&&?&\cr\tablerule
&&1&&(1,3,5,7)&&$z_0^{15},z_1^5,z_2^3,z_3^2z_0,\ldots
(19)$&&15&&9&&?&\cr\tablerule
&&1&&(1,3,5,8)&&$z_0^{16},z_1^5z_0,z_2^3z_0,z_3^2,\ldots
(20)$&&16&&10&&?&\cr\tablerule
&&1&&(2,3,5,9)&&$z_0^{9},z_1^6,z_2^3z_1,z_3^2,\ldots
(13)$&&18&&7&&Y&\cr\tablerule &&1&&(3,3,5,5)&&$g_{\ss{5}}(z_0,z_1),
f_{\ss{3}}(z_2,z_3)$&&15&&5&&Y&\cr\tablerule
&&1&&(3,5,7,11)&&$z_0^{6}z_2,z_1^5,z_2^2z_3,z_3^2z_0,\ldots
(8)$&&25&&5&&Y&\cr\tablerule
&&1&&(3,5,7,14)&&$z_0^{7}z_2,z_1^5z_0,g_{\ss{2}}(z_2^2,z_3),\ldots
(9)$&&28&&6&&Y&\cr \tablerule
&&1&&(3,5,11,18)&&$g_{\ss{2}}(z_0^6,z_3),z_1^5z_2,z_2^3z_0,\ldots
(10)$&&36&&6&&Y&\cr \tablerule
&&1&&(5,14,17,21)&&$z_0^{7}z_3,z_1^4,z_2^3z_0,z_3^2z_1,z_0^5z_1z_2$&&56&&4&&Y&
\cr\tablerule &&1&&(5,19,27,31)&&$z_0^{10}z_3,z_1^4z_0,z_2^3,z_3^2z_1,
z_0^7z_1z_2$&&81&&3&&Y&\cr\tablerule
&&1&&(5,19,27,50)&&$z_0^{20},z_0^{10}z_3,z_3^2,z_1^5z_0,z_2^3z_1,
z_0^7z_1^2z_3$&&100&&4&&Y&\cr \tablerule
&&1&&(7,11,27,37)&&$z_0^{10}z_1,z_1^4z_3,z_2^3,z_3^3z_0$&&81&&3&&Y&\cr
\tablerule
&&1&&(7,11,27,44)&&$z_0^{11}z_1,z_2^3z_0,z_1^8,z_1^4z_3,z_3^2,z_0^4z_1^3z_2$&&88&&4&&Y&\cr \tablerule &&1&&(9,15,17,20)&&$z_0^{5}z_1,z_1^4,z_2^3z_0,z_3^3$&&60&&3&&Y&\cr\tablerule &&1&&(9,15,23,23)&&$z_0^{6}z_1,z_1^4z_0,z_2^3, z_2^2z_3,z_3,z_2z_3^2,z_3^3$&&69&&5&&Y&\cr \tablerule
&&1&&(11,29,39,49)&&$z_0^{8}z_2,z_1^4z_0,z_2^3,z_3^2z_1$&&127&&3&&Y&\cr
\tablerule
&&1&&(11,49,69,128)&&$z_0^{17}z_2,z_1^5z_0,z_2^4,z_2^2z_3,z_3^2$&&256&&2&&Y&\cr
\tablerule
&&1&&(13,23,35,57)&&$z_0^{8}z_1,z_1^4z_2,z_2^2z_3,z_3^2z_0$&&127&&3&&Y&\cr
\tablerule
&&1&&(13,35,81,128)&&$z_0^{17}z_1,z_1^5z_2,z_2^3z_0,z_3^2$&&256&&2&&Y&\cr
\tablerule
&&2&&(2,3,4,5)&&$z_0^6,z_1^4,z_2^3,z_3^2z_0,\ldots
(10)$&&12&&5&&?&\cr\tablerule
&&2&&(2,3,4,7)&&$z_0^7,z_1^4z_0,z_2^3z_0,z_3^2,\ldots (11)$&&14&&6&&?&
\cr\tablerule &&2&&(3,4,5,10)&&$z_0^5z_2,z_1^5,z_2^4,z_3^2,\ldots
(9)$&&20&&5&&Y&\cr\tablerule
&&2&&(3,4,6,7)&&$g_{\ss{3}}(z_0^2,z_2),z_1^3z_2,z_3^2z_1,z_0z_3z_1^2,z_0^2z_1^3$&&18&&6&&?& 
\cr \tablerule
&&2&&(3,4,10,15)&&$z_0^{10},z_1^5z_3,z_2^3,z_3^2,\ldots(10)$&&30&&7&&Y&\cr\tablerule
&&2&&(3,7,8,13)&&$\!\!\!\!\!z_0^7z_2,z_1^3z_2,z_2^2z_3,z_3^2z_0,z_0^5z_1,
z_0^3z_1z_3, z_0^2z_1z_2^2\!\!\!\!\!$&&29&&5&&?&\cr\tablerule
&&2&&(3,10,11,19)&&$\!\!\!\!\!\!z_0^{10}z_3,z_1^3z_2,z_2^2z_3,z_3^2z_0,
z_0^7z_1^2,z_0^4z_1z_3,z_0^3z_1z_2^3\!\!\!\!\!\!$&&41&&5&&?&\cr \tablerule
&&2&&(5,13,19,22)&&$z_0^{7}z_3,z_1^4z_0,z_2^3,z_3^2z_1,z_0^5z_1z_2$&&57&&3&&Y&\cr\tablerule &&2&&(5,13,19,35)&&$z_0^{14},z_0^7z_3,z_3^2,z_1^5z_0,z_2^3z_1,z_0^5z_1^2z_2$&&70&&3&&Y&\cr
\tablerule
&&2&&(6,9,10,13)&&$z_0^{6},z_1^4,z_2^3z_0,z_3^2z_2,z_0^3z_1^2$&&36&&4&&Y&\cr
\tablerule
&&2&&(7,8,19,25)&&$z_0^{7}z_1,z_1^4z_3,z_2^3,z_3^2z_0,z_0^2z_1^3z_2$&&57&&3&&Y&\cr\tablerule
&&2&&(7,8,19,32)&&$z_0^{8}z_1,z_1^8,z_1^4z_3,z_3^2,z_2^3z_0,z_0z_2^3,
z_0^3z_1^3z_2$&&64&&4&&Y&\cr
\tablerule
&&2&&(9,12,13,16)&&$z_0^{4}z_1,z_1^4,z_2^3z_0,z_3^3$&&48&&3&&Y&\cr\tablerule
&&2&&(9,12,19,19)&&$z_0^{5}z_1,z_1^4z_0,z_2^3,z_2^2z_3,z_2z_3^3,z_3^3$&&57&&5&&Y&\cr
\tablerule
&&2&&(9,19,24,31)&&$z_0^{9},z_1^3z_2,z_2^3z_0,z_3^2z_1$&& 81&&3&&Y&\cr
\tablerule
&&2&&(10,19,35,43)&&$z_0^{7}z_2,z_1^5z_0,z_2^3,z_3^2z_1$&&105&&3&&Y&\cr
\tablerule
&&2&&(11,21,28,47)&&$z_0^{7}z_2,z_1^5,z_2^3z_1,z_3^2z_0$&&105&&3&&Y&\cr
\tablerule
&&2&&(11,25,32,41)&&$z_0^{6}z_3,z_1^3z_2,z_2^3z_0,z_3^2z_1$&&107&&3&&Y&\cr
\tablerule
&&2&&(11,25,34,43)&&$z_0^{10},z_1^4z_0,z_2^2z_3,z_3^2z_1$&&111&&3&&Y&\cr
\tablerule
&&2&&(11,43,61,113)&&$z_0^{15}z_2,z_1^5z_0,z_2^3z_1,z_3^2$&&226&&2&&Y&\cr
\tablerule
&&2&&(13,18,45,61)&&$z_0^{9}z_1,z_1^5z_2,z_2^3,z_3^2z_0$&&135&&3&&Y&\cr
\tablerule}}
}
\vfil\eject
\centerline{
\vbox{\tabskip=0pt \offinterlineskip
\def\tablerule{\noalign{\hrule}}
\halign to380pt {\strut#& \vrule#\tabskip=1em plus2em&
     \hfil#& \vrule#& \hfil#& \vrule#&  \hfil#& \vrule#&
     \hfil#& \vrule#&
     \hfil#& \vrule#&
     \hfil#& \vrule#\tabskip=0pt\cr\tablerule
\omit&height2pt&\multispan{11}&\cr
&&\multispan{11}\hfil {\bf Table 1.} (cont.)
Sporadic Examples of $\calz_{\bf w}$ of Index $1\leq I\leq10$\hfil&\cr\tablerule
&&\omit\hidewidth Index\hidewidth&&
\omit\hidewidth ${\bf w}$\hidewidth&&
\omit\hidewidth Monomials of $f_{\bf w}$\hidewidth&&
\omit\hidewidth $d$\hidewidth&&
\omit\hidewidth $b_2$\hidewidth&&
\omit\hidewidth K-E\hidewidth&\cr\tablerule
&&2&&(13,20,29,47)&&$z_0^{6}z_3,z_1^3z_3,z_2^3z_1,z_3^2z_0$&&107&&3&&Y&\cr
\tablerule
&&2&&(13,20,31,49)&&$z_0^{7}z_1,z_1^4z_2,z_2^3z_3,z_3^2z_0$&&111&&3&&Y&\cr
\tablerule
&&2&&(13,31,71,113)&&$z_0^{15}z_1,z_1^5z_2,z_2^3z_0,z_3^2$&&226&&2&&Y&\cr
\tablerule
&&2&&(14,17,29,41)&&$z_0^{5}z_3,z_1^5z_0,z_2^2z_3,z_3^2z_1$&&99&&3&&Y&\cr
\tablerule
&&3&&(5,7,11,13)&&$z_0^4z_3,z_1^4z_0,z_2^3,z_3^2z_1,z_0^3z_1z_2$&& 
33&&3&&?&\cr\tablerule
&&3&&(5,7,11,20)&&$z_0^{8},z_1^5z_0,z_2^3z_1,z_3^2,z_0^4z_3,z_0^3z_1^2,z_2^2$&&40&&4&&Y&
\cr\tablerule
&&3&&(11,21,29,37)&&$z_0^{6}z_2,z_1^4z_0,z_2^2z_3,z_3^2z_1$&&95&&3&&Y&\cr
\tablerule
&&3&&(11,37,53,98)&&$z_0^{13}z_2,z_1^5z_0,z_2^3z_1,z_3^2$&&196&&2&&Y&\cr
\tablerule &&3&&(13,17,27,41)&&$z_0^{6}z_1,z_1^4z_2,z_2^2z_3,z_3^2z_0$&&
95&&3&&Y&\cr \tablerule
&&3&&(13,27,61,98)&&$z_0^{13}z_1,z_1^5z_2,z_2^3z_0,z_3^2$&&196&&2&&Y&\cr
\tablerule
&&3&&(15,19,43,74)&&$z_0^{7}z_2,z_1^7z_0,z_2^3z_1,z_3^2$&&148&&2&&Y&\cr
\tablerule
&&4&&(5,6,8,9)&&$z_0^{3}z_3,z_1^4,z_2^3,z_3^2z_1,z_0^2z_1z_2$&&24&&3&&?&\cr\tablerule
&&4&&(5,6,8,15)&&$z_0^{6},z_1^5,z_2^3z_1,z_3^2,z_0^3z_3,z_0^3z_3,z_0^2z_1^2z_2$&&30&&4&&?&
\cr\tablerule
&&4&&(9,11,12,17)&&$z_0^{5},z_1^3z_2,z_2^3z_0,z_3^2z_1$&&45&&3&&?&\cr\tablerule
&&4&&(10,13,25,31)&&$z_0^{5}z_2,z_1^5z_0,z_2^3,z_3^2z_1$&&75&&3&&Y&\cr\tablerule
&&4&&(11,17,20,27)&&$z_0^{4}z_3,z_1^3z_2,z_2^3z_0,z_3^2z_1$&&71&&3&&?&\cr
\tablerule
&&4&&(11,17,24,31)&&$z_0^{5}z_2,z_1^4z_0,z_2^2z_3,z_3^2z_1$&&79&&3&&Y&\cr
\tablerule
&&4&&(11,31,45,83)&&$z_0^{11}z_2,z_1^5z_0,z_2^3z_1,z_3^2$&&166&&2&&Y&\cr
\tablerule
&&4&&(13,14,19,29)&&$z_0^{4}z_2,z_1^3z_3,z_2^3z_1,z_3^2z_0$&&71&&2&&?&\cr
\tablerule
&&4&&(13,14,23,33)&&$z_0^{5}z_1,z_1^4z_2,z_2^2z_3,z_3^2z_0$&&79&&3&&Y&\cr
\tablerule
&&4&&(13,23,51,83)&&$z_0^{11}z_1,z_1^5z_2,z_2^3z_0,z_3^2$&&166&&2&&Y&\cr
\tablerule
&& 5&&(11,13,19,25)&&$z_0^4z_2,z_1^4z_0,z_2^2z_3,z_3^2z_1$&&63&&3&&?&\cr
\tablerule
&&5&&(11,25,37,68)&&$z_0^{9}z_2,z_1^5z_0,z_2^3z_1,z_3^2$&&136&&2&&Y&\cr
\tablerule
&&5&&(13,19,41,68)&&$z_0^{9}z_1,z_1^5z_2,z_2^3z_0,z_3^2$&&136&&2&&Y&\cr
\tablerule
&&6&&(7,10,15,19)&&$z_0^{5}z_1,z_1^3z_2,z_2^3,z_3^2z_0$&&45&&3&&?&\cr\tablerule
&&6&&(11,19,29,53)&&$z_0^{7}z_2,z_1^5z_0,z_2^4,z_3^2$&&106&&2&&Y&\cr\tablerule
&&6&&(13,15,31,53)&&$z_0^{7}z_1,z_1^5z_2,z_2^3z_0,z_3^2$&&106&&2&&Y&\cr
\tablerule
&&7&&(11,13,21,38)&&$z_0^{5}z_2,z_1^5z_0,z_2^3z_1,z_3^2$&&76&&2&&Y&\cr\tablerule
&&8&&(7,11,13,23)&&$z_0^{5}z_1,z_1^3z_2,z_2^3z_0,z_3^2$&&46&&2&&?&\cr\tablerule
&&8&&(7,18,27,37)&&$z_0^{9}z_1,z_1^3z_2,z_2^3,z_3^2z_0$&&81&&3&&?&\cr\tablerule
&&9&&(7,15,19,32)&&$z_0^{7}z_1,z_1^3z_2,z_2^3z_0,z_3^2$&&64&&2&&?&\cr\tablerule
&&10&&(7,19,25,41)&&$z_0^{9}z_1,z_1^3z_2,z_2^3z_0,z_3^2$&&82&&2&&?&\cr\tablerule
&&10&&(7,26,39,55)&&$z_0^{13}z_1,z_1^3z_2,z_2^3,z_3^2z_0$&&117&&3&&?&\cr
\tablerule}} } \vskip -4pt \eightrm \ \ \ \ \ \ \ \ \ \ \ \ \ \ \ \ \ \ \ \ \
\  * (for lack of space only the total number of monomial terms  in $f_\bfw$
is indicated) \tenrm
\bigskip

The computer program indicates that there are neither series solutions 
nor sporadic solutions satisfying the 
hypothesis of Theorem \alg.5 for $I>10$. In fact, an easy argument shows that
there are no such solutions for sufficiently large $I$.
\footnote{*}{\ninerm The code for the C program used to generate the tables of
the Theorem \alg.5 are available at the following URL
http://www.math.unm.edu/$\tilde{\phantom{o}}$galicki/papers/publications.html.
} \bigskip

\bigskip
\noindent{\sc Remarks} \alg.6: As already mentioned
when $2I=w_0+w_1$ or $2I<3w_0$ the sufficient condition for the existence of
the K\"ahler-Einstein metric fails. In all such cases
techniques of [DK] say nothing about the corresponding $\calz_w$.
Table 2 below gives the classical examples of $\calz_w$ that can be
written as smooth hypersurfaces in weighted projective spaces and which by
other methods admit K-E metrics [TY]. We only indicate a particular polynomial
as the most general such hypersurface involves many such monomials which for
reasons of space we do not include. We are not aware of any non-smooth orbifold
examples for which the conditions of Theorem \alg.5 fail, but are known to
admit K\"ahler-Einstein metrics by another method.

\bigskip
\centerline{
\vbox{\tabskip=0pt \offinterlineskip
\def\tablerule{\noalign{\hrule}}
\halign to380pt {\strut#& \vrule#\tabskip=1em plus2em&
     \hfil#& \vrule#& \hfil#& \vrule#&  \hfil#& \vrule#&
     \hfil#& \vrule#&
     \hfil#& \vrule#\tabskip=0pt\cr\tablerule
\omit&height2pt&\multispan{9}&\cr
&&\multispan{9}\hfil {\bf Table 2.}
Classical del Pezzo surfaces as 
hypersurfaces in $\bbc\bbp^3(\bfw)$\hfil&\cr\tablerule
&&\omit\hidewidth Index\hidewidth&&
\omit\hidewidth ${\bf w}$\hidewidth&&
\omit\hidewidth $f_{\bf w}$\hidewidth&&
\omit\hidewidth $d$\hidewidth&&
\omit\hidewidth $\calz_\bfw$\hidewidth&\cr\tablerule
&&1&&(1,1,1,1)&&$z_0^3+z_1^3+z_2^3+z_3^3$
&&3&&$\bbc\bbp(2)\#6\overline{\bbc\bbp(2)}\phantom{{\matrix{1\cr2\cr}}}$&\cr\tablerule
&&1&&(1,1,1,2)&&$z_0^4+z_1^4+z_2^4+z_3^2$
&&4&&$\bbc\bbp(2)\#7\overline{\bbc\bbp(2)}\phantom{{\matrix{1\cr2\cr}}}$&\cr\tablerule
&&1&&(1,1,2,3)&&$z_0^6+z_1^6+z_2^3+z_3^2$
&&6&&$\bbc\bbp(2)\#8\overline{\bbc\bbp(2)}\phantom{{\matrix{1\cr2\cr}}}$&\cr\tablerule
&&2&&(1,1,1,1)&&$z_0^2+z_1^2+z_2^2+z_3^2$
&&2&&$\bbc\bbp(1)\times\bbc\bbp(1)\phantom{{\matrix{1\cr2\cr}}}$&\cr\tablerule
&&3&&(1,1,1,1)&&$z_0+z_1+z_2+z_3$
&&1&&$\bbc\bbp(2)\phantom{oooo{\matrix{1\cr2\cr}}}$&\cr\tablerule}}
}

\bigskip
\baselineskip = 10 truept
\centerline{\bf \pr.  Proof of Theorem \alg.5}  
\bigskip

The proof goes as follows: 
we begin by proving two lemmas which by Lemma \alg.2 will reduce the problem
to one in which the integers $m_2$ and $m_1$ of \alg.1 are bounded. The
computer programs then give a printout of solutions of \alg.1 which give all
log del Pezzo orbifolds $\calz_\bfw$
such that $(\calz_\bfw ,D)$ might be klt for all  
$D \equiv -{2+\gre \over 3} K_{\calz_\bfw}$.  
In the case
that $m_0$ is unbounded, these are the {\it series} solutions, and the
case when $m_0$ is bounded, they are the {\it sporadic} solutions.
Corollary \link.7 is sufficient, for most of the sporatic solutions, to check
that $(\calz_\bfw ,D)$ is klt for all $D \equiv -{2+\gre \over 3}
 K_{\calz_\bfw}$.  
However, for the series solutions the klt condition is
much more tedious. Indeed, in several cases we are unable to give an answer.
The results are collected in tables. In the last column of the tables we
indicate whether the klt condition holds (Y) or is unknown (?). Since there
exists a K\"ahler-Einstein orbifold metric whenever
$(\calz_\bfw,D)$ is klt for all $D \equiv -{2+\gre \over 3}
 K_{\calz_\bfw}$, this column is labeled K-E. The second to
the last column gives the second Betti number of $\calz_\bfw$ 
which is
computed by the method described in section \top. Also, where space permits,
we give the monomials that make up the defining weighted homogeneous
polynomial. This is indicated in the tables for most of the sporadic
solutions. 

\noindent{\sc Lemma} \pr.1: \tensl  Suppose $2I \geq 3w_0$.  Then  
for any $\gre > 0$, there exists $D \in |-K_{\calz_\bfw}|$ such
that $(\calz_\bfw,{2+\gre \over 3}D)$ is not klt. \tenrm

\noindent{\sc Proof}: Consider the section $z_0 \in H^0(\calz_\bfw,
\calo_{\calz_\bfw}(w_0))$ 
and let $D$
be its zero divisor.  Then by hypothesis $rD \in |\calo_{\calz_\bfw}
(-K_{\calz_\bfw})|$ for $r
\geq {3\over 2}$.   But then ${2+\epsilon\over  3}D$ can never be klt at a
generic point $y\in D$, for $\epsilon > 0$, 
 since it has multiplicity greater than 1. \hfill\za

\noindent{\sc Lemma} \pr.2: \tensl Suppose $\calz_\bfw 
 \subset \bbp(w_0,w_1,w_2,w_3)$ 
is a hypersurface of index $I$ with $2I=w_0+w_1.$
Then for any $\gre > 0$, there exists $D \in |-K_{\calz_\bfw}|$ such
that $(\calz_\bfw,{2+\gre \over 3}D)$ is not klt. \tenrm

\noindent{\sc Proof}: First we notice in the proof of Lemma \alg.2 that
$2I=w_0+w_1$ occurs in precisely one case and in this case we have
$w_3=w_1+w_2-I.$ But then $\bfw$ must have the form 
$$\bfw =(I-n,I+n,w,w+n)$$
for some $w\in \bbz^+$ with $w\geq I+n,$ and some non-negative integer $n<
I$ and in this case the degree $d=2w+n+I.$ Now suppose that
$z_0^{a_0}z_1^{a_1}z_2^{a_2}z_3^{a_3}$ is a monomial occurring with nonzero
coefficient in the polynomial defining $\calz_\bfw$. 
We claim that if $a_0 = 0$ then
$a_1 \neq 0$.  Indeed, if $a_0 = a_1 = 0$ then
$$ a_2w + a_3w = 2w + 2I + n. $$
But $ w \geq I + n$ so the only possible solutions would require $a_2 + a_3 =
3$ and $w = 2I + n$ and one readily checks that these hypersurfaces are
never quasi-smooth.
Thus the divisor $D = \{z_0 = 0\} \cap \calz_\bfw $ 
has at least two components, $E$ and
$F,$ where $E$ is the line $z_0 = z_1 = 0$ and $F$ is defined, inside the
weighted projective plane $\{z_0 = 0\}$, by a polynomial $f(z_1,z_2,z_3)$.  
Moreover, $f(z_1,z_2,z_3) = z_2^2 + z_1^2(g(z_1,z_2,z_3))$.  
Note that the point $P = (0,0,0,1) \in \calz_\bfw$ 
since $w+n$ does not divide the
degree of $\calz_\bfw $.  
Thus if 
$\pi: \bbc^2 \ra{1.3} \calz_\bfw $ is a local cover of the quotient singularity
at $P = (0,0,0,1)$  then $\pi^\ast D = \pi^\ast E + \pi^\ast F$ has 
multiplicity at least $1 + {\rm mult}_0(\pi^\ast F)$ at the origin.  
To compute the multiplicity of $\pi^\ast F$ at the origin, let $Y = Z(f(z_1,
z_2,z_3))$.  Then $Y \cap \calz_\bfw  = F \cup G,$ 
where $G$ does not contain the point 
$(0,0,0,1)$.    Thus 
$$
{\rm mult}_0(\pi^\ast F) = {\rm mult}_0(\pi^\ast Z(f)) \geq 2.
$$
Consequently,
$(\calz_\bfw,{2 + \gre \over 3}D)$ is never klt as it always has multiplicity
$> 2$ at $0$.  But ${I\over I-n}D \in |-K_{\calz_\bfw}|$ so
this completes the proof of the lemma. \hfill\za

The analysis of most of the sporadic examples of Table 1 is easily done
with help of Corollary \link.7 which can restated for this purpose
as:

\noindent{\sc Corollary} \pr.3: \tensl Let 
$\bfw=(w_0,w_1,w_2,w_3)$ and $\calz_\bfw\subset\bbp(\bfw)$ be a quasi-smooth
surface of degree $d=w_0+w_1+w_2+w_3-I$. Then $\calz_\bfw$ admits
a K\"ahler-Einstein metric if $2Id<3w_0w_1$. If the line
$(z_0=z_1=0)\not\subset\calz_\bfw$ then $2Id<3w_0w_2$ is also
sufficient. If the point $(0,0,0,1)\not\in \calz_\bfw$ then $2Id<3w_0w_3$ is also
sufficient.
\tenrm

We begin with the more complicated analysis of the
infinite series examples.  In what follows $\bbp$ will frequently
denote the appropriate
weighted projective space, with the weights $\bfw$ being understood; similarly
$\calz_\bfw \subset \bbp$ will denote a hypersurface in the weighted 
projective space, depending on an integer parameter $k$ for all of the
series examples.

\bigskip
\noindent$\bullet$
We consider now the hypersurface $\calz_\bfw\subset
{\bbp}(4,2k+1,2k+1,4k)$ 
given by the zero set of the homogeneous
polynomial of degree $d=8k+4$
$$
f_\bfw(z_0,z_1,z_2,z_3)= z_0^{2k+1} + z_0z_3^2 + g(z_1,z_2),
$$
where $g$ is of degree $4$ and we assume that the polynomial $g$
has 4 distinct  roots on the projective line.  Thus
$-K_{\calz_\bfw} = \calo_{\calz_\bfw}(2)$.  The
only singularities of $\calz_\bfw$ lie at the singular
points of ${\bbp}(\bfw)$.  
We can apply Corollary \link.7 to see that $(\calz_\bfw,D)$ is klt at all smooth
points of ${\bbp}(\bfw)$, for all $k>0$, for all $D 
\equiv -K_{\calz_\bfw}$ since
$$
d \cdot \kappa = (8k+4)2 < 4(2k+1)(2k+1).
$$

\medskip
We now pass to the singular points of $\calz_\bfw$ which are slightly more involved
as we need to consider a desingularization.  We first treat the point
$P=(0,0,0,1)$.  We would like to apply the preceding discussion and use 
Shokurov's inversion of adjunction (Lemma \est.5)
for which we need to verify that 
${\rm mult}_0(\pi^\ast D) \leq 2$
for any effective $\bbq$--divisor $D \equiv -K_{\calz_\bfw}$: 
here $\pi: (\bbc^2,0)
\rightarrow 
(\calz_\bfw,P)$ is a local cover of the quotient singularity at $P$.
We will consider the linear
series $|z_1,z_2|$ on $\calz_\bfw$.  This has only isolated base points on $\calz_\bfw$,
including the point $P$ in question, and
consequently if $E$ is a general member of this system then $D \cap E$ will
be proper.  Pulling back to $\bbc^2$ to compute the intersection multiplicity
at $0$ we find, since the ramification of $\pi$ over $P$ is of degree $4k$,
$$
{\rm mult}_0(\pi^\ast D) \leq 4k D \cdot E
= {2(2k+1)(8k+4)4k\over 4(2k+1)^2 4k} = 2,
$$ 
and thus we may apply Lemma \est.5.  Hence
it is sufficient, in order to establish that
$(\calz_\bfw,D)$ is log--canonical at $P$, to show
that
$$
p_\ast^{-1}(\pi^\ast D)|E 
$$
is a sum of points each with coefficient at most one.  

Consider the
hypersurface $z_0 = 0$ on $\calz_\bfw$.  This is a union of four lines,
$L_1,L_2,L_3,L_4$, one for each (projective)
zero of the polynomial $g(z_1,z_2)$.  Thus we have
$$L_1+L_2+L_3+L_4 \in |\calo_{\calz_\bfw}(4)|.\leqno{\pr.4}$$
Note that if an effective divisor $D \equiv -K_{\calz_\bfw}$ does not 
contain any of the lines $L_i$, then computing $(L_1+L_2+L_3+L_4)D$
will immediately establish that $(\calz_\bfw,D)$ is klt at $P$.  Thus the 
representatives of $-K_{\calz_\bfw}$ which are of concern are those
containing some or all of the lines $L_i$.  So suppose
$$D \equiv a_1L_1 + a_2L_2 + a_3L_3 + a_4L_4 + D^\prime,\leqno{\pr.5}$$
where $D^\prime$ meets each of the four lines properly.  We are interested in
bounding the $a_i$ which can be accomplished intersection theoretically.  In
particular we compute
$$L_i \cdot L_j = {1\over 4k}\qquad{\rm for}\qquad i \not= j.
\leqno{\pr.6}$$
To see this, note that the four lines $L_i$ are algebraically equivalent.
Thus, using \pr.4 above gives

$$\calz_\bfw \cdot \calo_{{\bbp}}(4)  \cdot \calo_{\bbp}(2k+1) = 
4(L_i \cdot \calo_{\bbp}(2k+1))\qquad
{\rm for \ all}\  i.\leqno{\pr.7}$$
The left hand side of \pr.7 is ${4(8k+4)(2k+1)\over 4(2k+1)^24k} =
{1\over k}$.  On the other hand one can check, choosing the appropriate
representative for $\calo_{\bbp}(2k+1)$ that $\calo_{\bbp}(2k+1) \cdot 
\calz_\bfw  =
L_j + C,$ where $C$ does not meet $L_i$;
\pr.6 follows immediately.

Intersecting \pr.4 with $\calo_{\bbp}(1)$ gives
$$\calo_{{\bbp}}(1) \cdot L_i = {1\over 4(2k+1)k}, \qquad \forall \ i.
\leqno{\pr.8}$$

Using \pr.6 and \pr.8 we can compute $L_i^2$:
$$L_i^2 = L_i \cdot \calo_{\bbp}(4) - \sum_{j \neq i} L_i \cdot L_j  = 
    {1-6k\over (2k+1)4k}.\leqno{\pr.9}$$
We have, by \pr.8
$$D \cdot L_i = {1\over 2k(2k+1)}.\leqno{\pr.10}$$
Using \pr.5 we obtain from \pr.10
$$a_i L_i^2 + \sum_{j \neq i} a_j L_j \cdot L_i
+ D^\prime \cdot L_i
= {1\over 2k(2k+1)}.\leqno{\pr.11}$$ 
To compute the terms in \pr.11 observe first that $\sum_{i=1}^4 a_i \leq 2$;
this follows from the fact that ${\rm mult}_0(\pi^\ast D) \leq 2$.  
Next, note that 
$$
D^\prime \cdot L_i \leq
D \cdot D = \calo_{\calz_\bfw}(2) \cdot \calo_{\calz_\bfw}(2) = {4(8k+4)\over
4(2k+1)^24k} = {1\over k(2k+1)}.
$$
Finally, using \pr.6 and \pr.9, \pr.11 becomes
$$
{a_i(1-6k)\over (2k+1)4k} + {1\over 4k}\sum_{j \neq i}a_j + {1\over k(2k+1)}
\geq {1\over 2k(2k+1)}.
$$
Clearing denominators then yields
$$
a_i \leq {k+1\over 2k}.
$$
Thus if $k \geq 2$ we have $a_i \leq {3\over 4}$. 

\medskip
Next we bound ${\rm mult}_0(\pi^\ast D^\prime)$.  Since $L_i \cap D^\prime$ is
proper for any $i$ we compute, using \pr.8
$$
{\rm mult}_0(\pi^\ast D^\prime) \leq 4k L_i \cdot D \leq 
{2\over 2k+1}.
$$
Returning to the inversion of adjunction set--up, we see that 
$p_\ast^{-1}(\pi^\ast(D))$ is a sum of four points
(corresponding to the pull--backs of the four
of the lines $L_i$) having total multiplicity at most
${3\over 4}$ and another divisor with total degree at most ${2\over 2k+1}$.
A simple computation then shows that ${11\over 13}D$  is  klt at
$P = (0,0,0,1)$ for any $D \equiv -K_{\calz_\bfw}$.

\medskip
Next we turn to the points $P_i = (0,a_i,b_i,0)$ on $\calz_\bfw$; there are exactly 
four of these by our assumption on the homogeneous polynomial $f(z_1,z_2)$ 
and the singularities have index $2k+1$.
Since ${\rm BS}(\calo_{\calz_\bfw}(w_0w_1w_2) \otimes \cali_{P_i})$ 
is a finite set of points,
we will apply Lemma \est.12 with $A = \calo_{\bbp}(8k+4)$.  Then we have
$$
A \cdot D = {(8k+4)(2)(8k+4)\over 4(2k+1)(2k+1)4k}
          = {2\over k}.
$$
We can take $a = 4$ in Lemma \est.12 since the linear series 
$\calo_{\bbp}(8k+4)$
allows for an isolated zero at $P_i$ of multiplicity $4$.  
Thus for $k \geq 2$ we see that $(\calz_\bfw,{3\over4}D)$ is klt at $P_i$.  

\medskip
Finally, we deal with the points $Q_j = (a_j,0,0,b_j),$ where the quotient
singularities are of index $4$.  We again apply Lemma \est.12, this time
taking $A = \calo_{\bbp}(2k+1)$.   
This linear series has sections with an isolated
singularity of multiplicity $1$ at any of the $Q_j$.  We compute
$$
A \cdot D = {2(2k+1)(8k+4)\over 4(2k+1)^24k} = {1\over 2k}
$$
Thus for $k \geq 2$ we see that $(\calz_\bfw,{3\over 4}D)$ is
klt
at $Q_j$. 
Putting together 
all of our computations establishes that $(\calz_\bfw,{3\over 4}D)$ is
klt for $k \geq 2$ and for any $D \equiv -K_{\calz_\bfw}$.
\bigskip
\noindent$\bullet$
We consider now the hypersurface $\calz_\bfw \subset {\bbp}(3,3k+1,6k+1,9k+3)$
given by the zero set of the homogeneous
polynomial of degree $d=18k+6$
$$
f_\bfw(z_0,z_1,z_2,z_3)=z_0^{6k+2} + z_1^3z_3 +z_2^3z_0 +z_3^2.
$$
Here $-K_{\calz_\bfw}
 = \calo_{\calz_\bfw}(2)$.  The worst 
singularity in this case is at $x=(0,0,1,0),$ where the index is $6k+1$.  
Since ${\rm BS}(\calo_{\calz_\bfw}(w_0w_1) \otimes \cali_x)$
consists of finitely many points, Corollary \link.7 says that
$(\calz_\bfw,{5\over7}D)$ is klt at $x$ for $k \geq 1$ for any $D \equiv 
-K_{\calz_\bfw }$.

Next we consider the singular line $z_0 = z_2 = 0,$ where the index is $3k+1$.
This intersects $\calz_\bfw$ in the finitely many points.
Applying the third formula of Corollary \link.7 shows that
$(\calz_\bfw,{5\over7}D)$ is klt at each of these points for $k \geq 1$.

Finally, there are index $3$ singularities along the line $z_1 = z_2 = 0$
which again intersects 
$\calz_\bfw$ in finitely many points.  Applying Corollary
\link.7 to these 
shows that $(\calz_\bfw,{7\over8}D)$ is klt at these points for
$k \geq 1$.  Putting everything together we see that 
$(\calz_\bfw,{5 \over 7}D)$ 
is klt for any $D \equiv -K_{\calz_\bfw}$ provided $k \geq 1$.
\bigskip
\noindent$\bullet$
Next we consider the hypersurface $\calz_\bfw \subset
{\bbp}(6,6k+5,12k+8,18k+15)$
given by the zero set of the homogeneous
polynomial of degree $d=36k+30$
$$
f_\bfw(z_0,z_1,z_2,z_3)=z_0^{6k+5} + z_1^6 + z_1^3z_3 + z_0z_2^3 + z_3^2. 
$$
Here we have $\calo_{\calz_\bfw}(-K_{\calz_\bfw}) = 
\calo_{\calz_\bfw}(4)$.  This variety contains fewer
of the singular points of ${\bbp}$ than the previous example and
asymptotically in $k$,
the proportions of these weights are identical and
thus this also must be $\gra$-klt for all $k$ sufficiently large and for 
appropriate $\gra$.  More specifically,
we must deal with one point at infinity, $x = (0,0,1,0)$ which has index
$12k+8$.  Using $\calo_{\calz_\bfw}(6(6k+5))$ in Corollary \link.7 we see that
$(\calz_\bfw,{5 \over 7} D)$ is klt at $x$.

Next there is the line $z_0=z_2=0,$ where the index of ${\bbp}$ is $6k+5$.
Applying Corollary \link.7 as above gives $(\calz_\bfw,{5 \over 7}D)$ 
klt for  $k \geq 1$.  
Finally, there are is the line $z_1=z_2=0$ of index $3$.  
It follows once more from Corollary \link.7 that 
$(\calz_\bfw,{7 \over 8}D)$ is klt 
at each of these points for all $k \geq 2$.  There are singularities
of index 2 along the line $z_1=z_3=0$; $(\calz_\bfw,D)$ is 
klt at these
points for all $k\geq1$. Thus we find that for all
$k \geq 2$, $(\calz_\bfw,{5 \over 7}D)$ is klt for any $D \equiv
 -K_{\calz_\bfw}$.  

\bigskip
\noindent$\bullet$
Next we consider the hypersurface $\calz_\bfw
={\bbp}(3,3k+1,3k+2,3k+2)$
given by the zero set of the homogeneous
polynomial of degree $d=9k+6$
$$
f_\bfw(z_0,z_1,z_2,z_3)=z_0^{3k+2}+ z_0z_1^3 + g(z_2,z_3),
$$
where $g$ is a homogeneous polynomial of degree three 
with distinct zeroes.  The only singular points
to be analyzed in this case are
$(0,1,0,0)$ and $(0,0,a_i,b_i)$, where $(a_i,b_i)$ are the zeroes of $g$.
We first examine $P =(0,1,0,0)$, where ${\bbp}$
has a singularity of index $3k+1$.  Let $\pi: (\bbc^2,0) \rightarrow
(\calz_\bfw,P)$ be a local
cover  and
let $D \equiv -K_{\calz_\bfw}$.
Intersecting $D$ with a general 
member of $\calo_{\calz_\bfw}(3k+2)$ establishes that
$$ 
{\rm mult}_0(\pi^\ast D) \leq 2
$$
so that we can apply Shokurov's inversion of adjunction as before.  As before
we write $L_1 + L_2 + L_3$ for the zero scheme of $z_0$ on $\calz_\bfw$.   We compute
$$
\eqalign{ L_i \cdot L_j&= {1\over3k+1}, \qquad i  \not= j, \cr
\calo_{\bf P}(1) \cdot L_i&= {1\over(3k+1)(3k+2)},\cr
L_i^2&= {-1-6k\over(3k+1)(3k+2).}\cr}$$
For $D \equiv -K_{\calz_\bfw}$ 
effective we write $D = a_1L_1+a_2L_2+a_3L_3 + D^\prime$
as before, with $D^\prime$ meeting the three lines properly,
 and we find, expanding $D \cdot L_i$ and using the estimates 
$D \cdot D^\prime \leq D \cdot D$, $a_1+a_2+a_3 \leq 2$
$$
{(6k+1)a_i\over(3k+1)(3k+2)} +{a_i-2\over3k+1} \leq {2\over(3k+1)(3k
+2)}.
$$
Expanding yields $a_i \leq {6k+6\over9k+3}$.  Moreover, intersecting with
$z_0 = 0$ shows that 
$$
{\rm mult}_0(\pi^\ast D^\prime) \leq {2\over 3(3k+2)}.
$$
Thus, with notation as before, $p_\ast^{-1}\pi^\ast D$ is a sum of points with
multiplicity at most ${6k+6\over9k+3} + {2\over3(3k+2)}$ 
and inversion of
adjunction applies to show that $(\calz_\bfw,{3 \over 4}D)$ 
is klt at $P$ for all $k \geq 1$.

The other singular points of $\calz_\bfw$ like along the line $z_0 = z_1 = 0$ and
these 
have index $3k+2$.  Let  $Q_i = (0,0,a_i,b_i) \in \calz_\bfw$ be one of these
three points.  
The analysis of $D$ above for $P = (0,1,0,0)$ applies at these points as
well, unchanged though exactly one of the lines $L_i$ passes through $Q_i$.
In order to justify the use of Lemma \est.5 we need to show that 
${\rm mult}_0(\pi^\ast D) \leq 2$ for a local cover of ${\bbp}$ at $Q_i$.
This can be done with Lemma \est.12 using the linear series 
$\calo_{\calz_\bfw}(3(3k+1))$ allowing for a singularity of multiplicity $3$ at
$Q_i$.  Thus we conclude that
$(\calz_\bfw,{3 \over 4}D)$ is klt  for all $k \geq 1$ and for all $D \equiv
-K_{\calz_\bfw}$.

\bigskip
\noindent$\bullet$
Next we consider the hypersurface $\calz_\bfw\subset
{\bbp}(6,6k+3,6k+5,6k+5)$
given by the zero set of the homogeneous
polynomial of degree $d=18k+15$
$$
f_\bfw(z_0,z_1,z_2,z_3)=z_0z_1^3 + z_0^{2k+2}z_1 +
g(z_2,z_3),
$$
where again $g$ 
is homogeneous of degree 3 and $\calz_\bfw$.
This passes through
two points at infinity $(1,0,0,0)$ and $(0,1,0,0)$ and also intersects
the singular lines $z_0 = z_1 = 0$ and $z_2 = z_3 = 0$ in finitely many points.
The analysis of $(0,1,0,0)$ and $(0,0,a_i,b_i)$ is identical to the prior
case, yielding $(\calz_\bfw, {5\over7}D)$ is klt at these points for all $k \geq 1$.
For the
points along the line $z_2 = z_3 = 0$, other than $(0,1,0,0)$,the index 
does not depend on $k$, and hence
these are certainly klt for $k$ sufficiently large;
more specifically, $(\calz_\bfw,{5\over7}D)$ is also klt at these points for $k \geq
2$.  Hence $(\calz_\bfw,{5 \over 7}D)$ is klt 
for all $k \geq 2$ and for all $D \equiv - K_{\calz_\bfw}$. 

\bigskip
\noindent$\bullet$
Lastly we consider the hypersurface $\calz_\bfw\subset
{\bbp}(3,3k,3k+1,3k+1)$
given by the zero set of the homogeneous
polynomial of degree $d=9k+3$
$$
f_\bfw(z_0,z_1,z_2,z_3)= z_0^{3k+1} + z_0z_1^3
+g(z_2,z_3),
$$
where $g$ is homogeneous of degree $3$.  This contains only 
one singular point at infinity $(0,1,0,0)$ as well as the points $(0,0,a_i,
b_i)$ and the points $(c_i,d_i,0,0)$.  Routine computations exactly as above
show that $(\calz_\bfw,{3 \over 4}D)$ is klt for $k \geq 2$ and for all
$D \equiv -K_{\calz_\bfw}$.  

\medskip  
\noindent{\sc Remark} \pr.12:
Let us make a few comments for the series which we have not
analyzed.  For the two series of index $6$, the simple
intersection theoretic argument
which establishes that ${\rm mult}_0(\pi^\ast D) \leq 2$ on a local cover
of one of the singular points, where $D \equiv -K_{\calz_\bfw}$, fails.  The 
desired bound may well still hold but more detailed analysis of the other
points of intersection of $D$ and the appropriate divisor of Corollary 
\link.7 would be necessary to establish this.  In the other cases which
we have not analyzed, the singularities of
``bad'' divisor $D = \{z_0 = 0\} \cap \calz_\bfw$ require more subtle 
analysis because in these cases the remaining three weights are distinct
and hence it is more difficult to compute the contribution of $D$ to the
tangent cone at the origin of the appropriate local
cover; when two of the
remaining three weights are the same, one quickly reduces to ${\bbp}(1,1,k)$
which is singular only at $(0,0,1)$ and it is easy to check which tangent
directions $D$ gives on the resolution.  In principal, however, each
of these cases could be decided with more involved computation.
\bigskip

\bigskip
\baselineskip = 10 truept
\centerline{\bf \top. The Topology of the Link $L_f$} 
\bigskip

Recall the well-known construction of Milnor [Mil] for isolated hypersurface
singularities: There is a fibration of $(S^{2n+1}-L_f)\ra{1.3} S^1$ whose fiber
$F$ is an open manifold that is homotopy equivalent to a bouquet of n-spheres
$S^n\vee S^n\cdots \vee S^n.$ The {\it Milnor number} $\mu$ of $L_f$ is the
number of $S^n$'s in the bouquet. It is an invariant of the link which
can be calculated explicitly in terms of the degree $d$ and weights
$(w_0,\ldots,w_n)$ by the formula [MO]
$$\mu =\mu(L_f)=\prod_{i=0}^n\bigl({d\over w_i}-1\bigr).\leqno{\top.1}$$
The closure $\bar{F}$ of $F$ is a manifold with boundary that is homotopy
equivalent  to $F,$ and whose boundary is precisely the link $L_f.$  Then the
topology of $L_f$ is determined by the {\it monodromy map} induced by the
$S^1_\bfw$ action. Milnor and Orlik [MO] use these facts to give an algorithm
for computing the Betti number $b_{n-1}(L_f)$ from the characteristic
polynomial $\grD(t)$ of the monodromy map. The procedure is this. Associate to
any monic polynomial $f$ with roots $\gra_1,\ldots,\gra_k\in \bbc^*$ its
divisor  
$$\hbox{div}~f= <\gra_1>+\cdots+<\gra_k>\leqno{\top.2}$$
as an element of the integral ring $\bbz[\bbc^*]$ and let $\grL_n = \hbox{div}
(t^n-1).$  The `rational weights' used in [MO] are just ${d\over w_i},$
and are written in irreducible form, ${d\over w_i}={u_i\over v_i}$. The divisor
of the characteristic polynomial is then determined by
$$\hbox{div}~\grD(t)~= \prod_i({\grL_{u_i}\over v_i}-1)~= 1+\sum
a_j\grL_j,\leqno{\top.3.}$$
where $a_j\in \bbz$ and the second equality is obtained by using  the
relations  $\grL_a\grL_b=\break\gcd(a,b)\grL_{lcm(a,b)}.$ The second Betti
number of the link is then given by 
$$b_2(L_f)=1+\sum_ja_j.\leqno{\top.4}$$

Furthermore, the following proposition was proved in [BG3]:

\noindent{\sc Proposition} \top.5: \tensl Let $L_f$ be the link
of an isolated singularity defined by a weighted homogeneous polynomial $f$ in
four complex variables with weights $\bfw.$ Suppose further that the weights
$\bfw$ are well-formed, then $\hbox{Tor}~(H_2(L_f,\bbz))=0.$ \tenrm    

Now a well-known theorem of Smale [Sm] says that any simply connected compact
5-manifold which is spin, and whose second homology group is torsion free, is
diffeomorphic to $S^5\#l (S^2\times S^3)$ for some non-negative integer $l.$
Furthermore, it is known [BG2,Mor] that any simply connected \Se manifold is
spin. Combining this with the development above gives

\noindent{\sc Theorem} \top.6: \tensl Let $L_f$ be the link associated to
a well-formed weighted homogeneous polynomial $f$ in four complex variables.
Suppose also that $L_f$ is spin, in particular, if $L_f$ admits an \Se metric.
Then $L_f$ is diffeomorphic to  $S^5 \# l(S^2\times S^3)$, where
$l=b_2(L_f)=1+\sum_ja_j.$  \tenrm

\bigskip
\baselineskip = 10 truept
\centerline{\bf \mod. The Moduli Problem for \Se 5-Manifolds}  
\bigskip

In this section we discuss the moduli problem for \Se manifolds. It is not our
intention here to discuss the general moduli problem for \Se structures, but
rather to present some results regarding the 5 dimensional non-regular examples
described in the previous sections together with the known regular \Se
5-manifolds. In what follows by moduli space we shall mean certain sets of
objects (sections of vector bundles) modulo the action of the group of
diffeomorphisms that are diffeotopic to the identity.  In the case of
orbifolds diffeomorphism means diffeomorphism in the sense of orbifolds. Thus,
in the case of complex moduli any two representatives are deformation
equivalent.

We begin by briefly discussing the case complex structures on the del
Pezzo surfaces obtained by blowing up $\bbc\bbp^2$ at $l$ distinct
points for $0\leq l\leq 8.$ We shall always assume that no two of these points
lie on a line nor any three lie on a conic which we refer to as {\it in general
position}. The complex structures on 
$\calz_l=\bbc\bbp^2\#l\overline{\bbc\bbp}^2$
are then determined by the complex coordinates of the $l$ points modulo the
action of the complex automorphism group $\gG(\calz_l).$ Since we can fix
precisely four points in general position with the action of $GL(3,\bbc),$ the
moduli space of complex structures $\calm_l^\bbc$ on $\calz_l$ for 
$l=0,\ldots, 4$ is a single point, whereas for $l=5,6,7,8$ it is an open
connected manifold of complex dimension $2(l-4)$. 

It is known that when $l>2$, there exists [Siu,TY,Ti1] a unique [BM], up to
complex automorphism, K\"ahler-Einstein structure associated to each
complex structure. It follows that in this case ($l>2$) there is a 1-1
correspondence between the complex structures on $\calz_l$ 
and the homothety classes of positive K\"ahler-Einstein
metrics modulo complex automorphisms. This gives the identification of moduli
spaces 
$$\calm^\bbc_l\simeq \calm^{KE}_l, \quad \hbox{for $l>2$} \leqno{\mod.1}$$ 
where $\calm^{KE}_l $ denotes the moduli space of homothety classes of
K\"ahler-Einstein metrics on $\calz_l.$
Moreover, when $l=1,2$ the space $\calm^{KE}_l $ is well-known to be empty,
whereas, $\calm_l^\bbc$ is a single point space. Furthermore, for $\bbc\bbp^2$
and $\bbc\bbp^1\times\bbc\bbp^1$ both $\calm_l^\bbc$ and $\calm^{KE}_l $ are
single point spaces. When $4<l<9$, it is possible that for two inequivalent
complex structures $J_1,J_2\in\calm_l^\bbc$ the Einstein metrics $g_1$, $g_2$
solving the corresponding Monge-Ampere equations coincide. However, it follows
from a theorem of Pontecorvo [Pon] that up to complex conjugation this can
happen only when the metric is anti-self-dual, that is, the anti-self-dual
part $W_+$ of the Weyl conformal tensor must vanish. But it is well-known (cf.
[Boy]) that a compact anti-self-dual K\"ahler surface has zero scalar
curvature. So this cannot happen for del Pezzo surfaces. The precise result is
given below in Proposition \mod.13, and the argument works equally well for
compact K\"ahler orbifolds.

In the more general case of log del Pezzo surfaces $\calz$ that is complex
compact surfaces with positive first Chern class and at most quotient
singularities much less is known. Here we are dealing with K\"ahler orbifolds
and as discussed earlier the existence problem of finding K\"ahler-Einstein
orbifold metrics is still open. However, the Bando-Mabuchi uniqueness theorem
carries over to the orbifold case, so when a positive K\"ahler-Einstein
orbifold metric exists it is unique up to homothety and complex automorphisms.
Thus, generally, if $\calm^{KE}$ denotes the moduli space of homothety classes
of K\"ahler-Einstein metrics on $\calz,$ and $\calm^\bbc$ denotes the moduli
space of complex structures on $\calz,$ then there is an injective map
$$\calm^{KE}\ra{1.8} \calm^\bbc. \leqno{\mod.2}$$

Next we relate the moduli of K\"ahler-Einstein structures on $\calz$ to the
moduli of \Se structures on $\cals.$ First we give a result in general
dimension  from our previous work [BG1], and then specialize to dimension five.

\noindent{\sc Proposition} \mod.3: \tensl Two $2n+1$ dimensional rank one
Sasakian  structures $\cals=(\xi,\eta,\Phi,g)$ and
$\cals'=(\xi',\eta',\Phi',g')$ are equivalent if and only if their
corresponding spaces of leaves $(\calz,\gro,J,h)$ and $(\calz',\gro',J',h')$
are equivalent as K\"ahler orbifolds. Furthermore, $\cals=(\xi,\eta,\Phi,g)$ is
\Se if and only if $(\calz,\gro,J,h)$ is K\"ahler-Einstein with scalar
curvature $4n(n+1).$  \tenrm 

\noindent{\sc Remark} \mod.4:  Any Sasakian structure
$\cals=(\xi,\eta,\Phi,g)$ has a canonically equivalent structure, namely the
{\it conjugate} Sasakian structure  $\bar{\cals}=(-\xi,-\eta,-\Phi,g).$ This
corresponds on the space of leaves to the complex conjugate K\"ahler
structure, namely $(\calz,-\gro,-J,h).$ 

Notice that the correspondence in Proposition \mod.4 is between rank one \Se
structures  and  K\"ahler-Einstein structures with a fixed scalar curvature, or
equivalently  homothety classes of K\"ahler-Einstein structures. Now there are
two possible types of deformations of \Se structures, those that deform the
foliation and those that do not. The latter all lie in one of the families
$\gF(\xi)$ discussed in section \trans, and the former correspond to different
base orbifolds $\calz$ (assuming both the original and deformed \Se structures
are rank one).  We believe that \Se structures lying in different families
$\gF(\xi)$ correspond to distinct components of the moduli space of \Se
structures, but we do not prove this here.

Now let us fix some notation. Let $\cals_l=S^5\#l(S^2\times S^3)$ denote by
$\calm^{SE}_l$ the moduli space of \Se structures on $\cals_l,$ and let
$\calm_l^{reg,SE}$ and $\calm_l^{1,SE}$ denote the moduli space of regular
\Se structures, and rank one \Se structures on $\cals_l,$ respectively. Then
we have natural inclusions 
$$\calm_l^{reg,SE}\subset \calm_l^{1,SE}\subset \calm_l^{SE}.\leqno{\mod.5}$$  
Combining our discussion above with [FK,BG1,Ti1-Ti3] we obtain for the regular
case: 

\noindent{\sc Proposition} \mod.6: \tensl The following hold:
\item{(i)} $\calm^{reg,SE}_l$ is not empty if and only if 
$l=0,1,3,4,5,6,7,8$.
\item{(ii)} $\calm^{reg,SE}_l$ is a single point when $l=0,1,3,4$.
\item{(iii)} $\calm^{reg,SE}_l$ is a connected complex manifold
of dimension $2(l-4)$ when $l=5,6,7,8.$ Furthermore, up to conjugation there
is precisely one Sasakian structure sharing the same Einstein metric $g.$ 
\tenrm
\bigskip

To analyze the moduli problem in the non-regular case, we 
begin by describing the group of complex automorphisms
$\gG_\bfw$ of the weighted
projective 3-space $\bbp(\bfw)$. We shall assume that
$\bbp(\bfw)$ is well-formed. Recall that
$\bbp(\bfw)$ can be defined as a scheme ${\rm Proj}(S(\bfw))$,
where 
$$S(\bfw)=\bigoplus_dS^d(\bfw)=\bbc[z_0,z_1,z_2,z_3].$$
The ring of polynomials $\bbc[z_0,z_1,z_2,z_3]$ is graded with grading
defined by the weights $\bfw=(w_1,w_1,w_2,w_3)$. As a projective variety
we can embed $\bbp(\bfw)\subset \bbc\bbp^N$ and then the
group  $\gG_\bfw$ is a subgroup of $PGL(N,\bbc)$. $\bbp(\bfw)$ is
a toric variety and we can describe $\gG_\bfw$ explicitly as follows:
Let $\bfw=(w_0,w_1,w_2,w_3)$ be ordered as before.
We consider the group $G(\bfw)$ of automorphisms of the graded ring $S(\bfw)$
defined on generators by
$$\varphi_\bfw\pmatrix{z_0\cr z_1\cr z_2\cr z_3\cr}=
\pmatrix{f_0^{(w_0)}(z_0,z_1,z_2,z_3)\cr
f_1^{(w_1)}(z_0,z_1,z_2,z_3)\cr f_2^{(w_2)}(z_0,z_1,z_2,z_3)\cr 
f_3^{(w_3)}(z_0,z_1,z_2,z_3)\cr},\leqno{\mod.7}$$
where $f_i^{(w_i)}(z_0,z_1,z_2,z_3)$ is an arbitrary weighted homogeneous
polynomial of degree $w_i$ in $(z_0,z_1,z_2,z_3)$. This is  
a finite dimensional Lie group and
it is a subgroup of $GL(N,\bbc)$.
Projectivising, we get $\gG_\bfw=\bbp_\bbc(G(\bfw))$.

Note that when $\bfw=(1,1,1,1)$ then $G(\bfw)=GL(4,\bbc)$.  Other than this
case three weights are never the same if $\bbp(\bfw)$ is well-formed. If two
weights coincide then $G(\bfw)$ contains $GL(2,\bbc)$ as a subgroup. Finally,
when all weights are distinct we can write
$$\varphi_\bfw\pmatrix{z_0\cr z_1\cr z_2\cr z_3\cr}=
\pmatrix{a_0z_0\cr
a_1z_1+f_1^{(w_0)}(z_1)\cr a_2z_2+f_2^{(w_2)}(z_0,z_1)\cr 
a_3z_3+f_3^{(w_3)}(z_0,z_1,z_2)\cr},\leqno{\mod.8}$$
where $(a_0,a_1,a_2,a_3)\in(\bbc^*)^4$ and $f_i^{(w_i)},\ \ i=1,2,3$
are weighted homogeneous polynomials of degree $w_i$. The simplest situation
occurs when $f_1=f_2=f_3$ are forced to vanish. Then $\gG_\bfw=(\bbc^*)^3$ is
the smallest it can possibly be as $\bbp(\bfw)$ is toric. 
This is, in fact, common to
many examples of the log del Pezzo suraces of Table 1. More
precisely, we have 

\noindent{\sc Lemma} \mod.9: \tensl Let $w_0<w_1<w_2<w_3$. If
for each $i=1,2,3$ the weight $w_i$ is not a $\bbz^+$-linear combination
of smaller weights with non-negative integer coefficients than
$\gG_\bfw=(\bbc^*)^3$.\tenrm

Let $S_\bfw^d\subset S(\bfw)$ be the vector subspace spanned by all monomials
in $(z_0,z_1,z_2,z_3)$ of degree $d=|w|-I$, and  
let $\hat{S}^d(\bfw)\subset S^d(\bfw)$ denote subset all quasi-smooth
elements, i.e. those polynomials $f_\bfw\in S_\bfw^d$ such that conditions
\link.4 hold. Then we define $m_\bfw^d$ to be the dimension of the subspace
generated by $\hat{S}_\bfw^d.$ Now the automorphism group $G(\bfw)$ acts on
$S_\bfw^d$ leaving the subset $\hat{S}^d(\bfw)$ of quasi-smooth polynomials
invariant. Thus, for each log del Pezzo surface that is a solution of Theorem
\alg.5 we define the moduli space  
$$\calm_\bfw^d=\hat{S}_\bfw^d/G(\bfw)=\bbp(\hat{S}_\bfw^d)/\gG_\bfw,
\leqno{\mod.10}$$   
with $n_\bfw^d={\rm dim}(\calm_\bfw^d).$ Now there is an injective map
$$\calm_\bfw^d\ra{2.0} \calm^\bbc(\calz_\bfw),\leqno{\mod.11}$$
and by our results of sections \alg~and \pr, each element in $\calm_\bfw^d$
correponds to a unique homothety class of K\"ahler-Einstein metrics modulo
$\gG_\bfw$ and hence, to a unique \Se structure on the corresponding 5-manifold
$\cals_l$ modulo the group $\gG_\bfw$ acting as CR automorphisms. The results
with non-trivial moduli are collected in Table 3 below.

\noindent{\sc Example} \mod.12: As an illustration, we 
shall calculate the moduli spaces of complex structures of
the three classical del Pezzo surfaces of Table 2. We begin with
the cubic in $\bbp(1,1,1,1)$. The subset $\hat{S}_{(1,1,1,1)}^3\subset
\bbc^{20}$ is a dense open complex submanifold. In this case
$G(\bfw)=GL(4,\bbc)$ and the quotient $\calm_{(1,1,1,1)}^3$ is a
complex manifold of dimension 4. It is well-known that in this case
we have the identification $\calz_f\simeq\calz_6$ so that $\calm_{(1,1,1,1)}^3$
can be identified with $\calm_6^\bbc$.

\noindent
The second example is
a degree 4 surface in $\bbp(1,1,1,2)$. Here
$\hat{S}^4(1,1,1,2)\subset
\bbc^{22}$ since the general weighted polynomial $f(\bfz)$ of degree 4 can be
written as  the sum
$$f(\bfz)=g^{(4)}(z_0,z_1,z_2)+g^{(2)}(z_0,z_1,z_2)z_3+\lambda z_3^2.$$
The group $G((1,1,1,2))$ is defined by
$$\varphi_\bfw\pmatrix{z_0\cr z_1\cr z_2\cr z_3\cr}=
\pmatrix{\bba\pmatrix{z_0\cr z_1\cr z_2\cr}\cr \alpha z_3+\phi^{(2)}(z_0,
z_1,z_2)\cr},\ \ \ \ \ \bba\in GL(3,\bbc),\ \ \alpha\in\bbc^*,$$
where $\phi^{(2)}(z_0,z_1,z_2)$ is an arbitrary homogeneous
polynomial of degree 2 in $(z_0,z_1,z_2).$ One can see that the
action is free and, hence, 
$\calm^4_{(1,1,1,2)}$ is a complex manifold of dimension
$n^4_{(1,1,1,2)}=6$. 
It is known that the smooth member $\calz_f\subset\bbp(1,1,1,2)$
is diffeomorphic to $\calz_7$ so that $\calm_{(1,1,1,2)}^4$
can be identified with $\calm_7^\bbc$.

\noindent
The third example is a degree 6
surface in $\bbp(1,1,2,3)$. Here
$\hat{S}^6(1,1,2,3)\subset
\bbc^{23}$ as the general weighted polynomial $f(\bfz)$ of degree 6 can be
written as the sum
$$f(\bfz)=g^{(6)}+g^{(4)}z_2+
g^{(2)}z_2^2+\lambda_1 z_2^3+g^{(3)}z_3+
g^{(1)}z_2z_3+\lambda_2z_3^2,$$
where $g^{(i)}=g^{(p)}(z_0,z_1)$ is a homogeneous polynomial
of degree $p=1,2,3,4,6$.
The group $G((1,1,2,3))$ in defined by
$$\varphi_\bfw\pmatrix{z_0\cr z_1\cr z_2\cr z_3\cr}=
\pmatrix{\bba\pmatrix{z_0\cr z_1\cr}\cr 
\alpha_2z_2+\phi_2^{(2)}(z_0,z_1)\cr \alpha_3z_3+\phi_3^{(3)}(z_0,
z_1)+\phi_3^{(1)}(z_0,z_1)z_2\cr},
\ \ \ \ \ \bba\in GL(2,\bbc),\ \ \alpha_2,\alpha_3\in\bbc^*,$$
and $\phi_i^{(p)}$, $i=1,2,3$ are homogeneous polynomials in
$(z_0,z_1)$ of degree $p$.
One can see that ${\rm dim}(G(1,1,2,3))=15$ and 
$\calm^6_{(1,1,2,3)}$ is a connected complex manifold of dimension
$n^6_{(1,1,2,3)}=8$. Again, the general smooth member
$\calz_f\subset\bbp(1,1,2,3)$ is diffeomorphic to $\calz_8$ so that
$\calm_{(1,1,2,3)}^4$ can be identified with $\calm_8^\bbc$. 

One can carry out similar calculations for any of the 
log del Pezzo surfaces. In Table 1 below we tabulate all the examples of log
del Pezzo surfaces of Theorem  \alg.5 which admit K\"ahler-Einstein metrics
and  for which $n_\bfw^d\geq1$ and indicate the corresponding link $L_f$ in
$S^7:$
\bigskip\bigskip
\centerline{
\vbox{\tabskip=0pt \offinterlineskip
\def\tablerule{\noalign{\hrule}}
\halign to380pt {\strut#& \vrule#\tabskip=1em plus2em&
     \hfil#& \vrule#& \hfil#& \vrule#&  \hfil#& \vrule#&
     \hfil#& \vrule#&
     \hfil#& \vrule#&
     \hfil#& \vrule#\tabskip=0pt\cr\tablerule
\omit&height2pt&\multispan{11}&\cr
&&\multispan{11}\hfil {\bf Table 3.}
Log del Pezzo surfaces with $n_\bfw^d\geq1$\hfil&\cr\tablerule
&&\omit\hidewidth Index\hidewidth&&
\omit\hidewidth ${\bf w}$\hidewidth&&
\omit\hidewidth $d$\hidewidth&&
\omit\hidewidth $m^d_\bfw$\hidewidth&&
\omit\hidewidth $n^d_\bfw$\hidewidth&&
\omit\hidewidth $\cals_l=L_f\subset S^7$\hidewidth&\cr\tablerule
&&1&&(2,2k+1,2k+1,4k+1)&&8k+4
&&12&&5&&$\#7(S^2\times S^3)$&\cr\tablerule
&&1&&(2,3,5,9)&&18
&&13&&5&&$\#6(S^2\times S^3)$&\cr\tablerule
&&2&&(3,4,10,15)&&30
&&10&&3&&$\#6(S^2\times S^3)$&\cr\tablerule
&&1&&(3,5,7,14)&&28
&&9&&4&&$\#5(S^2\times S^3)$&\cr\tablerule
&&1&&(3,5,11,18)&&36
&&10&&3&&$\#5(S^2\times S^3)$&\cr\tablerule
&&2&&(3,4,5,10)&&20
&&9&&3&&$\#4(S^2\times S^3)$&\cr\tablerule
&&1&&(3,5,7,11)&&25
&&8&&3&&$\#4(S^2\times S^3)$&\cr\tablerule
&&1&&(3,3,5,5)&&15
&&10&&2&&$\#4(S^2\times S^3)$&\cr\tablerule
&&2&&(7,8,19,32)&&64
&&7&&2&&$\#3(S^2\times S^3)$&\cr\tablerule
&&3&&(5,7,11,20)&&40
&&7&&2&&$\#3(S^2\times S^3)$&\cr\tablerule
&&1&&(5,14,17,21)&&56
&&5&&1&&$\#3(S^2\times S^3)$&\cr\tablerule
&&1&&(5,19,27,50)&&100
&&6&&1&&$\#3(S^2\times S^3)$&\cr\tablerule
&&1&&(7,11,27,44)&&88
&&6&&1&&$\#3(S^2\times S^3)$&\cr\tablerule
&&2&&(6,9,10,13)&&36
&&5&&1&&$\#3(S^2\times S^3)$&\cr\tablerule
&&1&&(5,19,27,31)&&81
&&5&&1&&$\#2(S^2\times S^3)$&\cr\tablerule
&&2&&(7,8,19,25)&&57
&&5&&1&&$\#2(S^2\times S^3)$&\cr\tablerule}}
}
\bigskip

Next we turn to the equivalence problem for the Einstein metrics. As mentioned
previously for the log del Pezzo surfaces with a Y in the last column of
Table 1 and the tables of Theorem \alg.5, there is a unique homothety class of
K\"ahler-Einstein metrics corresponding to each point of $\calm_\bfw^d.$ But
the question remains whether two inequivalent K\"ahler-Einstein structures can
share the same Riemannian metric. This is answered by a Theorem of Pontecorvo
[Pon] which can be restated for our purposes as:

\noindent{\sc Proposition} \mod.13: \tensl Let $(J,J')$ be two distinct log del
Pezzo structures on the same underlying orbifold $\calz$ that are both
compatible with the same Einstein metric $g.$ Then $J$ and $J'$ are complex
conjugates, i.e. $J'=-J.$ \tenrm   

\noindent{\sc Proof}: By Proposition 3.1 of [Pon] if $J'\neq -J$ then the
metric $g$ is anti-self-dual, i.e. $W_+=0.$ From [Boy] anti-self-dual
K\"ahler metrics must have vanishing scalar curvature, and this cannot happen 
for log del Pezzo surfaces which are positive. \hfill\za

Conjugate K\"ahler-Einstein structures correspond to the same point of the
moduli space $\calm_\bfw^d,$ so each point of $\calm_\bfw^d$ corresponds to 
a distinct diffeomorphism class of K\"ahler-Einstein metrics with scalar
curvature $4n(n+1).$ 
There are similar results on the Sasakian level  which are essentially due to
Tanno [Tan] and Tachibana and Yu [TaYu]. We reformulate this result as
follows:

\noindent{\sc Proposition} \mod.14: \tensl Let $M$ be a $(4n+1)$-dimensional
compact manifold. Let $\cals=(\xi,\eta,\Phi,g)$ and
$\cals'=(\xi',\eta',\Phi',g)$ be two distinct Sasakian structures on $M$
sharing the same Riemannian metric $g.$ Suppose further that $(M,g)$ is not a
sphere with the standard round Sasakian metric. Then $\cals'$ and $\cals$ are
conjugate Sasakian structures, i.e. $\cals'=(-\xi,-\eta,-\Phi,g).$ \tenrm

\noindent{\sc Proof}: Since $(M,g)$ is not a round sphere, a theorem of
Tachibana and Yu [TaYu] says that $g(\xi,\xi')$ is a constant, say $a.$
By Schwarz inequality $|a|\leq 1$ and since the Sasakian structures are
distinct we must have $-1\leq a< 1.$ If $a=-1$ then $\cals$ and $\cals'$ are
conjugate Sasakian structures. So assume $|a|<1,$ then following Tanno [Tan]
we define 
$$\xi^{''}= {\xi'-a\xi\over \sqrt{1-a^2}}.$$ 
Now we have 
$$g(\xi,\xi^{''})=0, \qquad g(\xi^{''},\xi^{''})=1,$$ 
and it follows from the characterization of Sasakian structures in terms of
the Riemannian curvature (cf. [BG2]) that $\xi^{''}$ defines a Sasakian
structure $\cals^{''}$ which is orthogonal to $\cals.$ But then this would
define a 3-Sasakian structure [Tan] which cannot exist in dimension $4n+1.$
Thus, $\cals$ and $\cals'$ are conjugate Sasakian structures. \hfill\za

Now combining Theorems \trans.23, \alg.5, \top.6, and Proposition \mod.14 with
the information in Tables 1 and 3 proves Theorem A of the Introduction.

\bigskip
\bigskip
\medskip
\centerline{\bf Bibliography}
\medskip
\font\ninesl=cmsl9
\font\bsc=cmcsc10 at 10truept
\parskip=1.5truept
\baselineskip=11truept
\ninerm
\item{[AFHS]} {\bsc B.S. Acharya, J. M. Figueroa-O'Farrill, C. M. Hull, C. M.,
and B. Spence}, {\ninesl Branes at conical singularities and holography},
Adv. Theor. Math. Phys. 2 (1998), 1249--1286.
\item{[B\"a]} {\bsc C. B\"ar}, {\ninesl Real Killing spinors and holonomy},
Comm. Math. Phys. 154 (1993), 509--52.
\item{[Be]} {\bsc A. Besse}, {\ninesl Einstein manifolds},
Springer-Verlag, Berlin-New York, 1987.
\item{[BFGK]} {\bsc H. Baum, T. Friedrich, R. Grunewald, I. Kath}, {\ninesl
Twistors and Killing spinors on Riemannian manifolds},
Teubner-Texte zur Mathematik, 124. 
B. G. Teubner Verlagsgesellschaft mbH, Stuttgart, 1991.
\item{[BG1]} {\bsc C. P. Boyer and  K. Galicki}, {\ninesl On Sasakian-Einstein
Geometry}, Int. J. of Math. 11 (2000), 873-909.
\item{[BG2]} {\bsc C. P. Boyer and  K. Galicki}, {\ninesl 3-Sasakian
Manifolds}, Surveys in Differential Geometry VI: 
{\it Essays on Einstein Manifolds};
A supplement to the Journal of Differential Geometry, pp. 123-184,
(eds. C. LeBrun, M. Wang); International Press, Cambridge (1999).
\item{[BG3]} {\bsc C. P. Boyer and  K. Galicki}, {\ninesl New Einstein Metrics
in Dimension Five},  math.DG/0003174, submitted for publication.
\item{[BM]} {\bsc S. Bando and T. Mabuchi}, {\ninesl Uniqueness of Einstein
K\"ahler Metrics Modulo Connected Group Actions}, Adv. Stud. Pure Math. 10
(1987), 11-40.
\item{[B\"o]} {\bsc C. B\"ohm}, {\ninesl~ Inhomogeneous Einstein metrics 
on low-dimensional spheres and other low-dimensional spaces},
Invent. Math. 134 (1998), 145-176.
\item{[Boy]} {\bsc C.P. Boyer}, {\ninesl Conformal Duality and Compact Complex
Surfaces}, Math. Ann. 274 (1986) 517-526.
\item{[DK]} {\bsc J.-P. Demailly and J. Koll\'ar}, {\ninesl Semi-continuity of
complex singularity exponents and K\"ahler-Einstein metrics on Fano
orbifolds}, preprint AG/9910118, to appear in Ann. Scient. Ec. Norm. Sup. Paris
\item{[Dol]} {\bsc I. Dolgachev}, {\ninesl Weighted projective varieties}, in
Proceedings, Group Actions and Vector Fields, Vancouver (1981) LNM 956, 34-71. 
\item{[ElK]} {\bsc A. El Kacimi-Alaoui}, {\ninesl Op\'erateurs transversalement
elliptiques sur un feuilletage riemannien et applications}, Compositio
Mathematica 79 (1990), 57-106. 
\item{[F]} {\bsc T. Friedrich}, {\ninesl 
Der erste Eigenwert des Dirac-Operators einer kompakten, 
Riemannschen Mannigfaltigkeit nichtnegativer Skalarkr\"ummung},
Math. Nachr. 97 (1980), 117-146.
\item{[FK]} {\bsc T. Friedrich and I. Kath}, {\ninesl Einstein manifolds of
dimension five with small first eigenvalue of the Dirac operator}, J. Diff.
Geom. 29 (1989), 263-279.
\item{[Fle]} {\bsc A.R. Fletcher}, {\ninesl Working with weighted complete
intersections}, Preprint MPI/89-95, revised version in  {\it Explicit
birational geometry of 3-folds},  A. Corti and M. Reid, eds., 
Cambridge Univ. Press, 2000,  pp 101-173.
\item{[JK1]} {\bsc J.M. Johnson and J. Koll\'ar}, {\ninesl K\"ahler-Einstein
metrics on log del Pezzo surfaces in weighted projective 3-space}, preprint
AG/0008129, to appear in Ann. Inst. Fourier.
\item{[JK2]} {\bsc J.M. Johnson and J. Koll\'ar}, 
{\ninesl Fano hypersurfaces in
weighted projective 4-spaces}, preprint AG/0008189, to appear in Experimental Math.
\item{[Ka]} {\bsc Y. Kawamata}, {\ninesl  The cone of curves of algebraic
varieties}, Annals of Math 119, 1984, pp. 603-33.
\item{[KMM]} {\bsc Y. Kawamata, K. Matsuda, and K. Matsuki}, {\ninesl
Introduction to the Minimal Model Problem}, Adv. Stud. Pure Math. 10 (1987),
283-360.
\item{[KM]} {\bsc J. Koll\'ar, and S. Mori}, {\ninesl Birational Geometry of
Algebraic Varieties}, Cambridge University Press, 1998.
\item{[Kol]} {\bsc J. Koll\'ar}, {\ninesl Rational Curves on Algebraic
Varieties}, Springer-Verlag, New York, 1996.
\item{[KW]} {\bsc I. R. Klebanov and E. Witten}, {\ninesl
Superconformal field theory on threebranes at a Calabi-Yau singularity},
Nuclear Phys. B 536 (1999), 199-218.
\item{[Ma]} {\bsc J. Maldacena}, {\ninesl
The large $N$ limit of superconformal field theories and supergravity},
Adv. Theor. Math. Phys. 2 (1998), 231-252.
\item{[Mil]} {\bsc J. Milnor}, {\ninesl Singular Points of Complex
Hypersurfaces}, Ann. of Math. Stud. 61, Princeton Univ. Press, 1968.
\item{[MO]} {\bsc J. Milnor and P. Orlik}, {\ninesl Isolated singularities
defined by weighted homogeneous polynomials}, Topology 9 (1970), 385-393. 
\item{[Mol]} {\bsc P. Molino}, {\ninesl Riemannian Foliations}, Progress in 
Mathematics 73, Birkh\"auser, Boston, 1988. 
\item{[Mor]} {\bsc S. Moroianu}, {\ninesl Parallel and Killing spinors on
$\hbox{Spin}^c$-manifolds}, Commun. Math. Phys. 187 (1997), 417-427.
\item{[MP]} {\bsc D. R. Morrison, M. R. Plesser}, {\ninesl
Non-spherical horizons I}, Adv. Theor. Math. Phys. 3 (1999), 1-81.
\item{[Na]} {\bsc A.M. Nadel}, {\ninesl Multiplier ideal sheaves
and existence of K\"ahler-Einstein metrics of positive scalar curvature}, Ann.
Math. 138 (1990), 549-596.
\item{[Pon]} {\bsc M. Pontecorvo}, {\ninesl Complex structures on Riemannian
four-manifolds}, Math. Ann. 309 (1997) 159-177.
\item{[Siu]} {\bsc Y.-T. Siu}, {\ninesl The existence of K\"ahler-Einstein
metrics on manifolds with positive anticanonical line bundle and a suitable
finite symmetry group}, Ann. Math. 127 (1988), 585-627.
\item{[Sm]} {\bsc S. Smale}, {\ninesl On the structure of 5-manifolds},
Ann. Math. 75 (1962), 38-46.
\item{[Tan]} {\bsc S. Tanno}, {\ninesl On the isometry groups of Sasakian
manifolds}, J. Math. Soc. Japan 22 (1970) 579-590.
\item{[TaYu]} {\bsc S. Tachibana and W.N. Yu}, {\ninesl On a Riemannian space
admitting more than one Sasakian structure}, T\^ohoku Math. J. 22
(1970) 536-540.
\item{[Ti1]} {\bsc G. Tian}, {\ninesl On Calabi's Conjecture for complex
surfaces with positive first Chern class}, Invent. Math. 101 (1990), 101-172.
\item{[Ti2]} {\bsc G. Tian}, {\ninesl K\"ahler-Einstein Manifolds with
positive scalar curvature}, Surveys in Differential Geometry VI:
{\it Essays on Einstein Manifolds};
A supplement to the Journal of Differential Geometry, pp.67-82,
(eds. C. LeBrun, M. Wang); International Press, Cambridge (1999).
\item{[Ti3]} {\bsc G. Tian}, {\ninesl On One of Calabi's Problems}, Proc.
Symp. Pure Math. 52 (1991) Part 2, 543-556.
\item{[Ton]} {\bsc Ph. Tondeur}, {\ninesl Geometry of Foliations}, Monographs
in Mathematics, Birkh\"auser, Boston, 1997. 
\item{[TY]} {\bsc G. Tian and S.-T. Yau}, {\ninesl K\"ahler-Einstein
metrics on complex surfaces with $c_1>0$}, Comm. Math. Phys. 112 (1987),
175-203.
\item{[Wi]} {\bsc E. Witten}, {\ninesl Anti
de Sitter space and holography}, Adv. Theor. Math. Phys. 2 (1998), 253-291.
\item{[Y]} {\bsc S. -T. Yau}, {\ninesl Einstein manifolds with zero Ricci
curvature},
Surveys in Differential Geometry VI:
{\it Essays on Einstein Manifolds};
A supplement to the Journal of Differential Geometry, pp.1-14,
(eds. C. LeBrun, M. Wang); International Press, Cambridge (1999).
\item{[YK]} {\bsc K. Yano and M. Kon}, {\ninesl
Structures on manifolds}, Series in Pure Mathematics 3, 
World Scientific Pub. Co., Singapore, 1984.
\medskip
\bigskip \line{ Department of Mathematics and Statistics
\hfil December 2000} \line{ University of New Mexico \hfil Revised, February
2001}  \line{ Albuquerque, NM 87131 \hfil } \line{ email: cboyer@math.unm.edu,
galicki@math.unm.edu, nakamaye@math.unm.edu\hfil} \line{ web pages:
http://www.math.unm.edu/$\tilde{\phantom{o}}$cboyer, 
http://www.math.unm.edu/$\tilde{\phantom{o}}$galicki \hfil}
\end